\documentclass{article}
\usepackage[english]{babel}

\usepackage[letterpaper,top=2cm,bottom=2cm,left=2cm,right=2cm,marginparwidth=1.75cm]{geometry}

\usepackage{calrsfs}
\usepackage{amssymb}
\usepackage{amsmath}
\usepackage{amsthm}
\usepackage{subcaption}
\usepackage{multirow}
\usepackage{graphicx}
\usepackage{amsfonts}
\usepackage{mathtools}
\usepackage{mathrsfs}
\usepackage{placeins}

\usepackage[mathscr]{euscript}
\usepackage[dvipsnames]{xcolor}
\usepackage[colorlinks=true, allcolors=OliveGreen]{hyperref}
\usepackage{filecontents}
\usepackage{tikz}
\usepackage{pgfplots}
\usetikzlibrary{external}
\usepackage{authblk}
\usepackage{todonotes}
\usepackage{array}
\newcolumntype{C}{>{\centering\arraybackslash}p{2.5cm}}
\newcommand{\averagel}{\{\!\!\{}
\newcommand{\averager}{\}\!\!\}}

\newcommand{\jumpl}{[\![}
\newcommand{\jumpr}{]\!]}

\newcommand{\partition}{\mathcal{T}_h}
\newcommand{\facesinternal}{\mathcal{F}^\mathrm{I}_h}
\newcommand{\faces}{\mathcal{F}_h}

\newcommand{\facesboundary}{\mathcal{F}^\mathrm{B}_h}

\newcommand{\Wh}{W_{h}^\mathrm{DG}}
\DeclareMathAlphabet{\mathcalligra}{T1}{calligra}{m}{n}

\mathtoolsset{showonlyrefs,showmanualtags}

\graphicspath{{Images/}}

\title{Exploring tau protein and amyloid-beta propagation: a sensitivity analysis of mathematical models based on biological data\footnote{\textbf{Funding}: Funded by the European Union (ERC SyG, NEMESIS, project number 101115663). Views and opinions expressed are, however, those of the author only and do not necessarily reflect those of the European Union or the European Research Council Executive Agency. Neither the European Union nor the granting authority can be held responsible for them. The author is a member of INdAM-GNCS. The present research is part of the activities of Dipartimento di Eccellenza 2023-2027.}}

\author[1]{Mattia Corti\footnote{mattia.corti@polimi.it}}

\affil[1]{MOX-Dipartimento di Matematica, Politecnico di Milano, Piazza Leonardo da Vinci 32, Milan, 20133, Italy}

\begin{document}
\maketitle

\begin{abstract}
Alzheimer's disease is the most common dementia worldwide. Its pathological development is well known to be connected with the accumulation of two toxic proteins: tau protein and amyloid-$\beta$. Mathematical models and numerical simulations can predict the spreading patterns of misfolded proteins in this context. However, the calibration of the model parameters plays a crucial role in the final solution. In this work, we perform a sensitivity analysis of heterodimer and Fisher-Kolmogorov models to evaluate the impact of the equilibrium values of protein concentration on the solution patterns. We adopt advanced numerical methods such as the IMEX-DG method to accurately describe the propagating fronts in the propagation phenomena in a polygonal mesh of sagittal patient-specific brain geometry derived from magnetic resonance images. We calibrate the model parameters using biological measurements in the brain cortex for the tau protein and the amyloid-$\beta$ in Alzheimer's patients and controls. Finally, using the sensitivity analysis results, we discuss the applicability of both models in the correct simulation of the spreading of the two proteins.
\end{abstract}
\section{Introduction}
\label{sec:introduction}
Neurodegenerative diseases represent a significant challenge in actual medical research. Due to the aging of the global population, the number of people affected by these pathologies is constantly increasing. Some of these pathologies are called proteinopathies due to the link between the disease progression and the aggregation and spreading of toxic proteins inside the central nervous system \cite{juckerSelfpropagationPathogenicProtein2013}. The principal example of proteinopathy is Alzheimer's disease, whose pathological development depends on the accumulation of two different proteins: tau protein and amyloid-$\beta$ \cite{scheltens_alzheimers_2021}. There is literature evidence that the alteration of these proteins' structure is connected to resistance to clearance mechanism, with consecutive agglomeration and neuronal death \cite{wilsonHallmarksNeurodegenerativeDiseases2023}. However, the delay between protein misfolding and the appearance of clinical symptoms is of the order of decades, complicating the construction of efficient medical treatments.
\par
First of all, Alzheimer's disease is associated with the accumulation of amyloid-$\beta$, which is considered the earliest hallmark of the disease. This protein is a standard product of Amyloid Precursor Protein (APP) metabolism, and it is also diffuse in the physiologic brain \cite{van_oostveen_imaging_2021}. However, its accumulation can lead to neuronal death. Secondly, Alzheimer's disease is also associated with neuronal and glial accumulation of misfolded tau protein. Tau is a microtubule protein of primary importance in the stabilization of neuronal cytoskeleton \cite{jouanne_tau_2017}. 
\par
Unless the importance of those proteins in the neurodegeneration process, their longitudinal monitoring is still problematic. Indeed, it can be performed only using positron emission tomography (PET) images, typically used only for diagnosis due to the procedure's high costs and invasive nature \cite{hampelCoreCandidateNeurochemical2008}. Other typical biomarkers studied in clinical practice are derived from chemical analysis of the cerebrospinal fluid (CSF). However, CSF studies require invasive procedures, and the correlation between concentrations of proteins inside the fluid and in the brain parenchyma still needs to be wholly understood \cite{hampelCoreCandidateNeurochemical2008}. In this context, mathematical models and numerical simulations can help in the longitudinal monitoring of protein propagation.
\par
Several mathematical models for prion dynamics have been proposed. Probably the most complete model, able to describe coagulation, fragmentation, and spreading processes, is the Smoluchowski model \cite{franchi-lorenzani}. However, the computational costs of solving this model have required the construction of more simplified ones, which try to maintain a good level of description, homogenizing the description of single proteins, oligomers, and plaques in a single term \cite{fornariPrionlikeSpreadingAlzheimer2019,weickenmeierPhysicsbasedModelExplains2019}. This work considers two mathematical models: the heterodimer and Fisher-Kolmogorov (FK) models. The heterodimer model \cite{matthaus_diffusion_2006} is based on a separate modeling of healthy and misfolded protein configurations. The model can explicitly represent proteins' production, destruction, and conversion. The FK model \cite{fisher-1937, kolmogorov-1937} is commonly used to describe the dynamics of biological systems. In this context, it provides a simplified description, valuable in slight variations in the healthy protein concentration \cite{weickenmeierPhysicsbasedModelExplains2019}. 
\par
In the protein spreading context, many different numerical methods are present in literature to describe the phenomenon, such as network diffusion models \cite{fornariPrionlikeSpreadingAlzheimer2019,thompsonProteinproteinInteractionsNeurodegenerative2020a}, finite element methods \cite{weickenmeierPhysicsbasedModelExplains2019}, and polytopal discontinuous Galerkin methods \cite{corti_discontinuous_2023,antonietti_heterodimer_2023,corti_positivity_preserving}. At the same time, many works studied the calibration of the model parameters starting from the final data, using inverse uncertainty quantification methodologies \cite{schaferBayesianPhysicsBasedModeling2021,schafer_correlating_2022,corti_uncertainty_2023}.
\par
This work aims to perform a sensitivity analysis of the models evaluating the impact of the variation of concentrations of the proteins in equilibrium conditions. This choice is driven by the importance of protein concentration levels in the medical literature. Indeed, it is deeply studied the patient-specificity of these biomarkers, to construct subgroups of patients \cite{duits_four_2021} and to control the impact of those values on disease development \cite{degerman_gunnarsson_high_2014}. Sensitivity analysis allows us to understand how the different input sources of uncertainty impact the uncertainty in the output of mathematical models. The values of protein concentrations are derived for the first time by constructing some probability distributions starting from biological \textit{post-mortem} measurements in the brain cortex \cite{wu_decrease_2011,roberts_biochemically-defined_2017}. In this case, a control on the quality of the resulting distribution has been visually performed by comparing the theoretical and experimental distribution functions. This procedure allows deriving specific parameter values for the tau protein and the amyloid-$\beta$ separately, highlighting the differences predicted by the models in the accumulation and spreading procedure for the different proteins. Moreover, performing the study on heterodimer and the FK model allows us to define the limits of applicability and understand when they could be considered equivalent and when they do not clearly describe the phenomenon. The results are compared with the medical literature, to evaluate the quality of the obtained longitudinal behaviors of total and phosphorylated tau, and of amyloid-$\beta$ and APP. Finally, the possible connections of the results with some oscillatory behaviors in cognitive assessments, CSF biomarkers \cite{mcdade_longitudinal_2018}, and PET measurement \cite{sanchez_longitudinal_2021,jagust_temporal_2021} have been discussed.
\par
\bigskip
The paper is organized as follows. Section~\ref{sec:mathematicalmodel} presents both models in their continuous strong formulation, and then Section~\ref{sec:stochasticparameters} highlights the stochastic parameter distributions. In Section~\ref{sec:methods}, we report the numerical methods adopted in the simulations. Finally, in Section~\ref{sec:results}, we report the results of the sensitivity analysis, we discuss them, and we report the limitations of the study in Section~\ref{sec:discussion}. Finally, in Section~\ref{sec:conclusions}, we draw some conclusions.
\section{Mathematical models in protheinopaties}
\label{sec:mathematicalmodel}
In this section, we briefly introduce the mathematical models to describe the processes of production, misfolding, clearance, and diffusion of the prionic proteins. In this work, we assume that the reaction parameters are dependent on a vector of parameters $\boldsymbol{v}\in \Gamma\subseteq\mathbb{R}^N$, where $\Gamma=\Gamma_1\times ...\times \Gamma_N$ is a hyperrectangle in which $\boldsymbol{v}$ is randomly varying. We specify later the exact form of the vector $\boldsymbol{v}$. For a final time $T>0$, the problem is dependent on time $t\in(0,T]$ and space $\boldsymbol{x}\in\Omega\subset\mathbb{R}^d$ ($d=2,3$). Here $\Omega$ is an open, bounded domain of $\mathbb{R}^d$. We denote by $(\cdot,\cdot)_\Omega$ the scalar product in $L^2(\Omega)$.
\subsection{The heterodimer model}
\label{subsec:2}
Firstly, we introduce the heterodimer model \cite{weickenmeierPhysicsbasedModelExplains2019}. The equations govern the dynamics of two different protein concentrations: the healthy $p= p(\boldsymbol{x},t,\boldsymbol{v})$ and misfolded $q = q(\boldsymbol{x},t,\boldsymbol{v})$ protein concentrations. This model describes the kinetics of prion pathogenesis incorporating prion conversion dynamics both considering monomeric and polymeric seeding hypotheses \cite{weickenmeierPhysicsbasedModelExplains2019}. Indeed, the system of equations reads:
\begin{equation}
\label{eq:hm_strong}
\begin{cases}
    \dfrac{\partial p}{\partial t}  =\nabla \cdot (\textbf{D} \nabla p) - k_{1} \, p - k_{12} \, p\, q + k_0 &  \mathrm{in} \: \Omega \times (0,T]\times\Gamma, \\[6pt]
    \dfrac{\partial q}{\partial t}  = \nabla \cdot (\textbf{D} \nabla q) - \tilde{k}_1\, q + k_{12}\, q \, p & \mathrm{in} \: \Omega \times (0,T]\times\Gamma, \\[6pt]
    (\textbf{D}\nabla p )\cdot \boldsymbol{n} = 0 \mathrm{,} \quad (\textbf{D} \nabla q )\cdot \boldsymbol{n}  = 0  & \mathrm{on} \; \partial\Omega\times (0,T]\times\Gamma, \\[6pt]
    p(\boldsymbol{x},0,\boldsymbol{v}) = p_{0}(\boldsymbol{x},\boldsymbol{v}) \mathrm{,} \quad q(\boldsymbol{x},0,\boldsymbol{v})  = q_{0}(\boldsymbol{x},\boldsymbol{v}) & \mathrm{in} \: \Omega\times\Gamma.
\end{cases}
\end{equation}
In Equation \ref{eq:hm_strong}, the first assumption we make is that the brain generates only healthy protein at a generation rate governed by the parameter $k_0=k_0(\boldsymbol{x},\boldsymbol{v})>0$. The biological clearance and destruction processes affect both healthy and misfolded proteins at rates $k_1=k_1(\boldsymbol{x},\boldsymbol{v})>0$ and $\tilde{k}_1=\tilde{k}_1(\boldsymbol{x},\boldsymbol{v})>0$, respectively. Finally, we have a conversion of proteins from healthy to misfolded state at a rate $k_{12}=k_{12}(\boldsymbol{x},\boldsymbol{v})>0$. To obtain a physically consistent description of the phenomenon, the reaction parameters we just introduced need to respect some relations \cite{thompsonProteinproteinInteractionsNeurodegenerative2020a}. Indeed, their values are directly related to the equilibrium values of the concentration in healthy and dementia conditions. Indeed, assuming a healthy subject, the concentration of misfolded protein should be $q_\mathrm{min}=0$, while the healthy ones should be around a value we call $p_\mathrm{max}=\frac{k_0}{k_1}$. These two values are the ones associated with the unstable equilibrium point of the heterodimer system $\boldsymbol{E}_1 = (p_\mathrm{max},0)$. On the contrary, the model also has a stable equilibrium point $\boldsymbol{E}_2 = (p_\mathrm{min},q_\mathrm{max})$, which is associated with an illness situation where we have reduced healthy protein concentration $p_\mathrm{min}=\frac{\tilde{k_1}}{k_{12}}$ and proliferation of misfolded ones $q_\mathrm{max}=\frac{k_0}{\tilde{k_1}}-\frac{k_1}{k_{12}}$.
\par
Concerning the diffusion part of the model, the spreading can be mathematically characterized by a tensor defined as follows \cite{weickenmeierPhysicsbasedModelExplains2019}:
\begin{equation}
    \label{eq:DiffusionTensor}
    \mathbf{D} = d_\mathrm{ext}\mathbf{I} + d_\mathrm{axn}\boldsymbol{\bar{a}} \otimes \boldsymbol{\bar{a}}.
\end{equation}
With this formulation, the first term models extracellular diffusion of magnitude $d_\mathrm{ext}=d_\mathrm{ext}(\boldsymbol{x})$, while the second term models anisotropic diffusion which happens along the axonal directions, denoted by unit vector $\boldsymbol{\bar{a}}=\boldsymbol{\bar{a}}(\boldsymbol{x})$ and of magnitude $d_\mathrm{axn}=d_\mathrm{axn}(\boldsymbol{x})$.
\par
Concerning the boundary conditions, we assume that parenchyma cannot exchange the proteins with the surrounding CSF, and so we impose homogeneous Neumann boundary conditions.
\subsubsection*{Bifurcations of heterodimer model equilibria}
\begin{figure}[t]
    \centering
    \includegraphics[width=\textwidth]{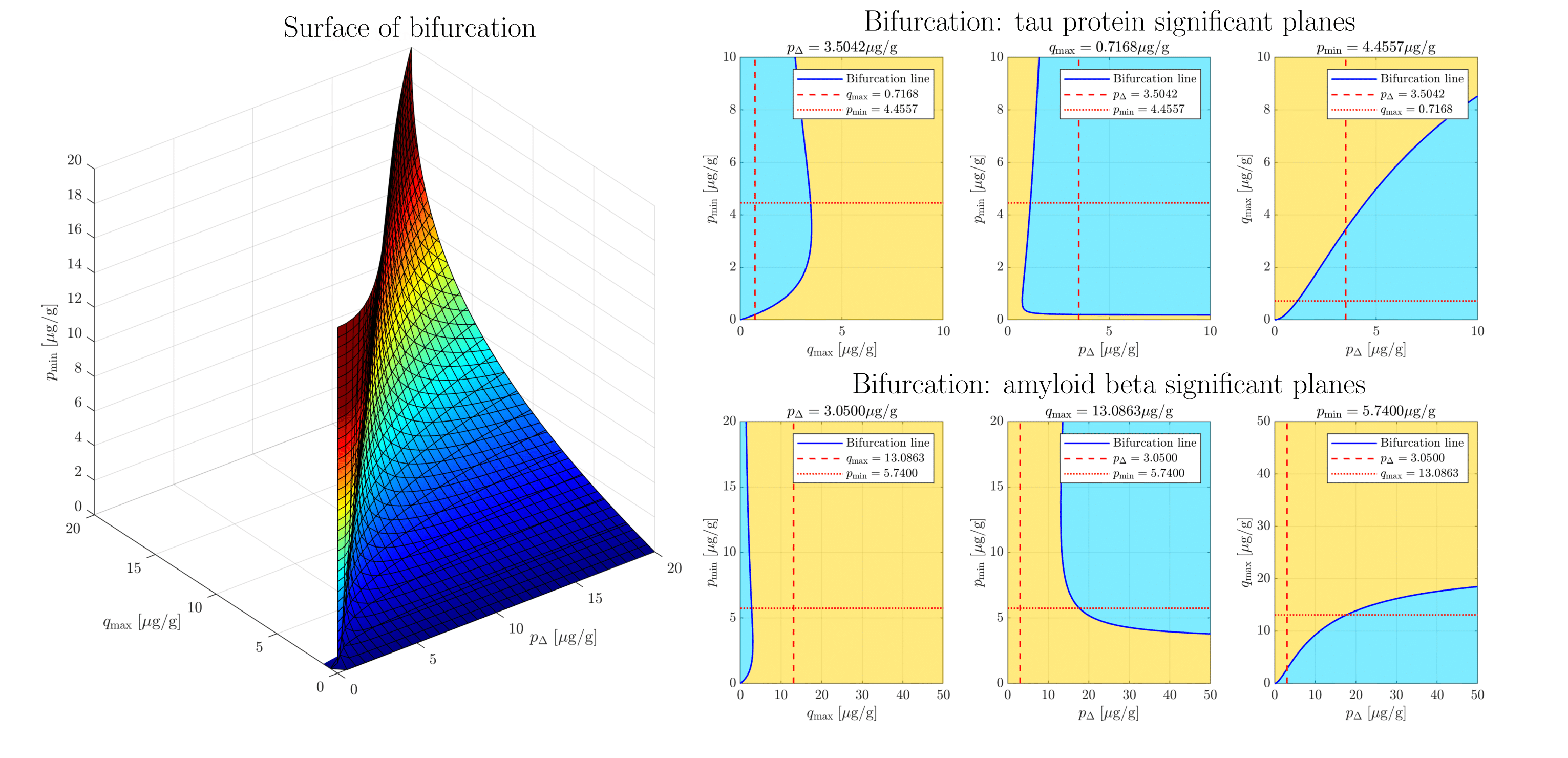}
   \caption{Bifurcation surface for the heterodimer model with respect to the parameters $p_\Delta$, $p_\mathrm{min}$, and $q_\mathrm{max}$ (left), and bifurcation lines intersected with the significant planes in the following simulations.}
   \label{fig:bifsurf}
\end{figure}

In the analysis of the equilibria of the heterodimer model, we can identify a bifurcation in the stable equilibrium $\boldsymbol{E}_2$. In particular, once we neglect the diffusion ($\mathbf{D}=\mathbf{0}$), the Jacobian associated with the problem \eqref{eq:hm_strong} evaluated at that point is:
\begin{equation}
    \mathbf{J}(\boldsymbol{E}_2) = 
    \begin{bmatrix}
        -q_\mathrm{max}\left(1-\dfrac{p_\mathrm{min}}{p_\Delta}\right) & -p_\mathrm{min} \\
        q_\mathrm{max} & 0
    \end{bmatrix},
\end{equation}
where $p_\Delta = p_\mathrm{max}-p_\mathrm{min}$. By computing the eigenvalues of the Jacobian matrix, we can observe that the real part of those is always negative provided the constraints defined above. However, depending on the parameters' values, they can be purely real (stable node), or complex and conjugate (stable focus) \cite{arnold2013dynamical}. The condition for this bifurcation is associated with the equation:
\begin{equation}
    (p_\Delta+p_\mathrm{min})^2 q_\mathrm{max} - 4 p_\mathrm{min} p_\Delta^2 = 0.
\end{equation}
This equation implicitly defines a bifurcation surface, which we report for the visualization in Figure \ref{fig:bifsurf}. Moreover, we report in the same figure, the bifurcation lines obtained by intersecting the surface with some planes that will be significant in the next simulations. In particular, in the image, the half-plane in light blue is associated with a stable focus behavior and the yellow part with a stable node one.
\subsection{The Fisher-Kolmogorov model}
A simplified formulation of problem~\eqref{eq:hm_strong} can be derived under certain assumptions. In particular, whenever the concentration of healthy proteins is much larger than the one of misfolded proteins, and if its variations can be considered negligible, we can perform a Taylor expansion of the function $p$ \cite{weickenmeierPhysicsbasedModelExplains2019, fornariPrionlikeSpreadingAlzheimer2019}. In this case, after defining a relative misfolded protein concentration as:
\begin{equation}
    c(\boldsymbol{x},t,\boldsymbol{v}) = \dfrac{q(\boldsymbol{x},t,\boldsymbol{v})}{q_{\mathrm{max}}},
\end{equation}
we can derive the FK model.
The variable $c=c(\boldsymbol{x},t,\boldsymbol{v})$ is a relative concentration assuming values in the interval $[0,1]$, where $0$ and $1$ means absence and high prevalence of misfolded proteins, respectively. Indeed, coherently with the starting model, $0$ is the unstable equilibrium, and $1$ is the stable one. Finally, the governing equation reduces to the FK model:
\begin{equation}
 \begin{cases}
     \dfrac{\partial c}{\partial t} =\nabla \cdot(\mathbf{D} \nabla\, c) + \alpha\,c(1-c),
    & \mathrm{in}\,\Omega\times(0,T]\times\Gamma,
    \\[8pt]
    (\mathbf{D}\nabla c) \cdot \boldsymbol{n} = 0, & 
    \mathrm{on}\;\partial \Omega \times(0,T]\times\Gamma,
    \\[8pt]
    c(\boldsymbol{x},0,\boldsymbol{v})=c_0(\boldsymbol{x},\boldsymbol{v}), & \mathrm{in}\;\Omega\times\Gamma.
    \\[8pt]
\end{cases}
\label{eq:fk_strong}
\end{equation}
where $\alpha = \alpha(\boldsymbol{x},\boldsymbol{v})$ is the conversion rate, describing misfolding, clearance, and protein production processes within a single term. Concerning the diffusion tensor $\mathbf{D}$, it is defined in Equation~\eqref{eq:DiffusionTensor}. The homogeneous Neumann boundary condition extends naturally to the concentration $c$ we introduce.
\section{Stochastic parameter distributions}
\label{sec:stochasticparameters}
A fundamental choice for performing the following investigations is the selection of the stochastic parameters that we will vary in the sensitivity analysis. In this work, we choose to analyze the impact of different values of maximum ($p_\mathrm{max}$,$q_\mathrm{max}$) and minimum ($p_\mathrm{min}$) concentrations of proteins on the final result. This choice is made mainly for two reasons: first, because these quantities are directly associated with the parameters of both models and second, because of the existence of clinical measures of them in literature. Moreover, these values are protein-dependent, allowing us to make some distinctions between the application of the models for different illnesses. For simplicity, we assume these parameters are constant in the domain $\Omega$. Concerning the concentration $p$, instead of using $p_\mathrm{max}$ as stochastic parameter, we use the difference $p_\Delta = p_\mathrm{max}-p_\mathrm{min}$. Indeed, this choice ensures that $p_\mathrm{max}>p_\mathrm{min}$ by simply selecting a positive probability distribution for the $p_\Delta$. Our modeling distribution is to choose for these parameters a Gamma distribution $\Gamma(a,b)$ such that the probability distribution is:
\begin{equation}
    \rho(y) = \dfrac{\beta^{\alpha+1}}{\Gamma(\alpha+1)}y^\alpha e^{-\beta y}.
\end{equation}
\par
Due to the presence of four reaction parameters in the equation \eqref{eq:hm_strong} and only three concentrations, we need a fixed parameter. Indeed, we evaluate the impact of different levels of protein concentrations on the dynamics at the same conversion rate.
\par
Finally, the components we choose for the definition of the vector are $\boldsymbol{v}=(p_\mathrm{min},p_\Delta,q_\mathrm{max})$. Introducing these quantities for the sensitivity analysis, we can rewrite the parameters of problems \eqref{eq:hm_strong} and \eqref{eq:fk_strong} as follows:
\begin{equation}
    k_0 = \dfrac{p_\mathrm{min}(p_\mathrm{min}+p_\Delta)q_\mathrm{max}}{p_\Delta}k_{12}, \qquad
    k_1 = \dfrac{p_\mathrm{min}q_\mathrm{max}}{p_\Delta}k_{12},
    \qquad
    \tilde{k}_1 = p_\mathrm{min}k_{12},
    \qquad
    \alpha = p_\Delta k_{12}.
\end{equation}
\par
The construction of parameter distributions needs to be protein-specific to adequately describe the physical and chemical processes connected to the pathophysiology of protein misfolding. Concerning the modeling choice in adopting a Gamma distribution for all the terms, this is mainly due to the necessity of constructing a positive distribution for the analyzed quantities. Moreover, we cannot fix a maximum value due to the nature of the measured quantities. For this reason, the choice of distribution with bounded support (such as the Beta distribution) is unsuitable for study purposes. Finally, the choice of a Gamma distribution is general because this distribution is a generalization of many others with those characteristics (i.e. exponential, or $\chi^2$-distribution) and does not have any assumption on the distribution of the logarithm, in contrast to a lognormal distribution.
\par
\subsection{Tau protein distributions}
Tau protein is a microtubule protein that plays a fundamental role in the stabilization of the neuronal cytoskeleton \cite{jouanne_tau_2017}. In Alzheimer's disease, this protein becomes hyperphosphorylated, with consequent disassociation from microtubules and aggregation into neurofibrillary tangles, one of the disease's most known hallmarks. This accumulation of phosphorylated tau proteins is associated with synaptic impairment and neurodegeneration \cite{drummond_phosphorylated_2020}.
\par
In this work, we estimate the parametrical distributions starting from some \textit{post-mortem} measurements of total and phosphorylated tau protein concentrations in the brain cortex for both Alzheimer's disease patients and controls \cite{wu_decrease_2011}. {In particular, starting from the raw data, we compute mean and variance and use them to derive the Gamma distribution parameters. We report the computed values in Table \ref{tab:tau_dist_param} with the specific reference of the figures from which we derived those values. We point out that to estimate the mean value $p_\Delta$, we derived the quantities from the control estimate of the total tau (which would represent $p_\mathrm{max}$). Then to obtain a meaningful value the mean of the $p_\mathrm{min}$ distribution is subtracted. For simplicity, we assume no variation in the variance estimate between $p_\mathrm{max}$ and $p_\Delta$.}
\par
\begin{figure}[t]
\centering
\includegraphics[width=\textwidth]{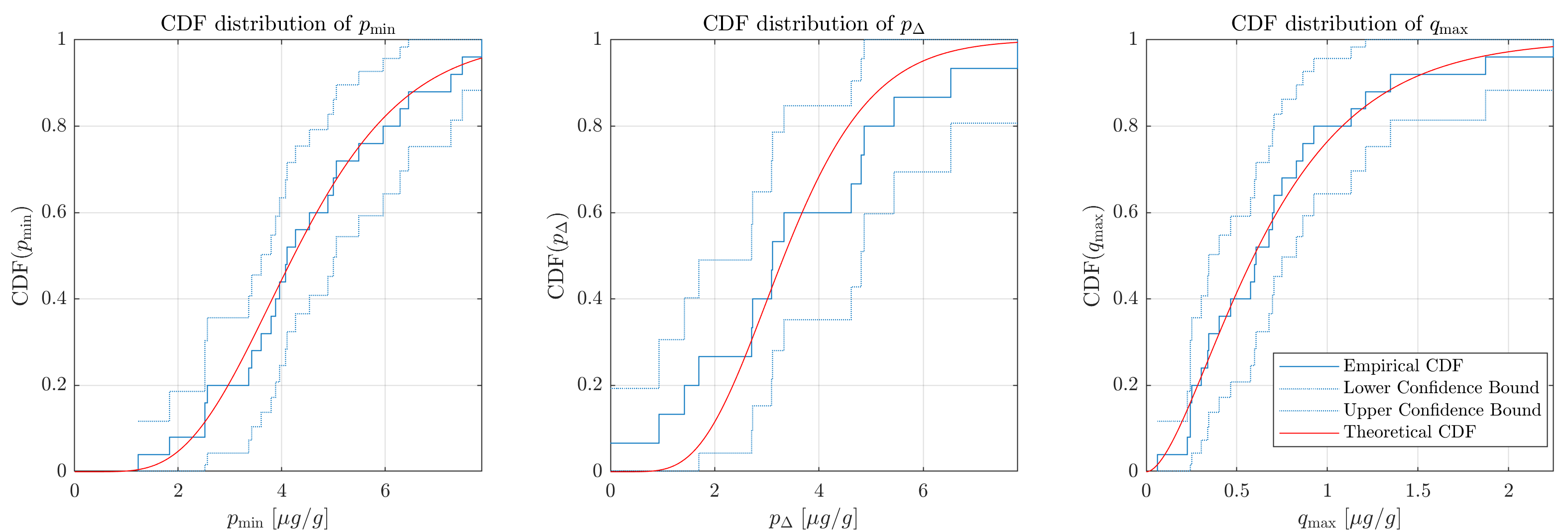}
\caption{{Comparison between empirical CDFs and theoretical CDFs for the model parameters in the tau protein simulations.}}
\label{fig:TauDistQual}
\end{figure}

{Finally, to confirm the quality of the theoretical distribution with respect to the raw data, we report the comparison between the Gamma cumulative distribution function (CDF) with respect to the empirical ones that can be derived from the raw data. As visible in Figure \ref{fig:TauDistQual}, the results confirm that the estimated distributions correctly fit the data and are contained in the $95\%$ confidence bounds.}
\begin{table}[t]
	\centering
	\begin{tabular}{|r|r l|r l|r l|r l|l|}
	\hline
	\multicolumn{1}{|c|}{\textbf{Parameter}} & \multicolumn{2}{c|}{\textbf{Mean}} & \multicolumn{2}{c|}{\textbf{Variance}} & \multicolumn{2}{c|}{\textbf{Parameter} $a$} & \multicolumn{2}{c|}{\textbf{Parameter} $b$} & \multicolumn{1}{c|}{{\textbf{Reference}}}\\ 
		\hline 
		 $p_\mathrm{min}\sim\Gamma(a,b)$  & $4.4557$ & $[\mu\mathrm{g/g}]$ & $3.0400$ & $[\mu\mathrm{g^2/g^2}]$ & $5.5307$ & $[-]$ & $1.4657$ & $[\mathrm{g}/\mu\mathrm{g}]$ & {\cite{wu_decrease_2011} Fig 5A (AD)} \\
		 $p_\Delta\sim\Gamma(a,b)$  & $3.5042$ & $[\mu\mathrm{g/g}]$ & $1.8217$ & $[\mu\mathrm{g^2/g^2}]$ & $5.7406$ & $[-]$ & $1.9236$ & $[\mathrm{g}/\mu\mathrm{g}]$& {\cite{wu_decrease_2011} Fig 5A (Control)} \\
		 $q_\mathrm{max}\sim\Gamma(a,b)$  & $0.7168$ & $[\mu\mathrm{g/g}]$ & $0.2737$ & $[\mu\mathrm{g^2/g^2}]$ & $0.8772$ & $[-]$ & $2.6189$ & $[\mathrm{g}/\mu\mathrm{g}]$ & {\cite{wu_decrease_2011} Fig 5B (AD)}  \\ \hline 
	\end{tabular}
	\caption{Estimated distributions parameters for the tau protein simulation \cite{wu_decrease_2011}.}
	\label{tab:tau_dist_param}
\end{table}
\subsection{Amyloid beta distributions}
Amyloid-$\beta$ is a protein produced through proteolytic processing of a transmembrane protein called APP. Amyloid-$\beta$ is produced when the APP is first cleaved by $\beta$-secretase. Suppose the produced protein undergoes a process of accumulation in the brain parenchyma. In that case, it can aggregate in plaques, whose formation is recognized as an early toxic event in the pathogenesis of Alzheimer's disease \cite{chen_amyloid_2017}. 
\par
\begin{figure}[t]
\centering
\includegraphics[width=\textwidth]{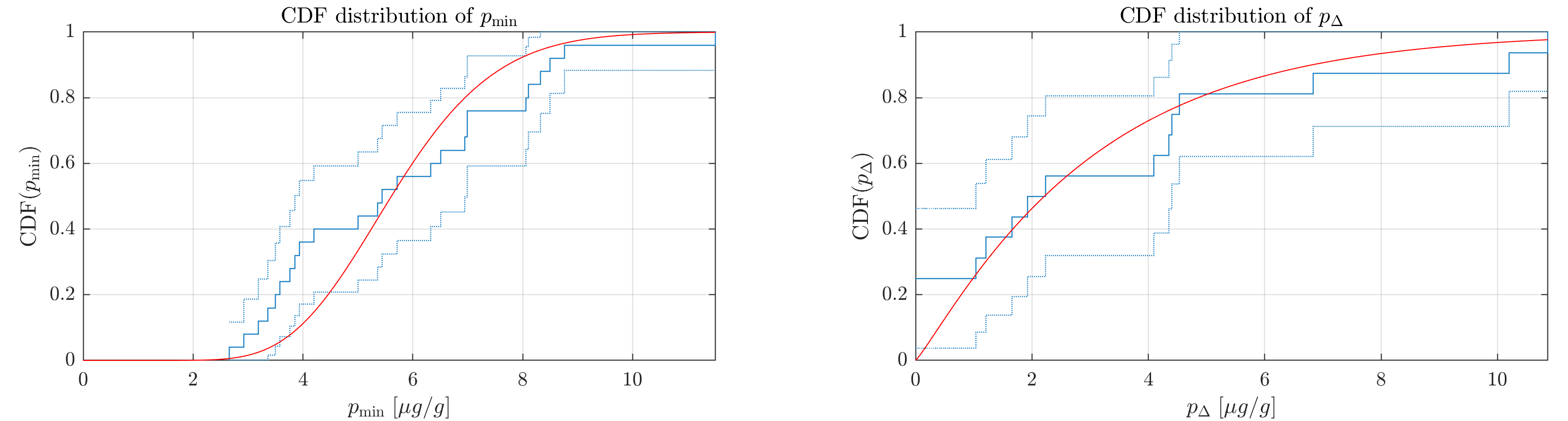}
\caption{{Comparison between empirical CDFs and theoretical CDFs for the model parameters in the amyloid-$\beta$ simulations.}}
\label{fig:BetaDistQual}
\end{figure}

In this work, we estimate the parametrical distribution of APP (in the model healthy protein $p$) starting from the \textit{post-mortem} cortex measurement in \cite{wu_decrease_2011}. {In particular, to reconstruct the distribution of $p_\mathrm{min}$ we consider the values of APP measurements reported in \cite{wu_decrease_2011} (Figure 4C) rescaling the values of the distribution in $\mu\mathrm{g/g}$ knowing that an approximation of the APP total mass can be $8.79 \mu\mathrm{g/g}$ obtained as the sum of the mean values of sAPP$\alpha$ and sAPP$\beta$ in the controls (considering a molecular weight of $100 \mathrm{kDa}$ as indicated in \cite{wu_decrease_2011}). The same procedure is applied considering the AD case, where the mean value estimated from sAPP$\alpha$ and sAPP$\beta$ is $5.74 \mu\mathrm{g/g}$. This approximation is chosen to consider both the products of $\alpha-$ and $\beta-$ secretase under which APP undergoes \cite{chen_amyloid_2017}. The obtained values are reported in Table \ref{tab:abeta_dist_param}. Also in this case, to estimate the mean value $p_\Delta$, we derived the quantities from the control estimate of the APP and then we subtracted the mean of the $p_\mathrm{min}$. In Figure \ref{fig:BetaDistQual}, we report the CDFs of the estimated distribution and the estimated CDFs from the data. We can observe that they correctly fit the data and are contained in the $95\%$ confidence bounds.}
\par
{The distribution associated with amyloid-$\beta$ is reconstructed starting from measurement in \cite{roberts_biochemically-defined_2017}. In this case, the values of mean and variance are derived from \cite{roberts_biochemically-defined_2017} (Table 1), where the mean and variance of total amyloid mass are divided by the average grey matter volume to obtain the values in Table \ref{tab:abeta_dist_param}. The choice of the Gamma distribution is done as an extension to what was done with the other quantities, however, a complete validation with the estimated CDFs is not possible due to the absence of raw data about the total mass of amyloid-$\beta$.}
\begin{table}[t]
	\centering
	\begin{tabular}{|r|r l|r l|r l|r l|l|}
	\hline
	\multicolumn{1}{|c|}{\textbf{Parameter}} & \multicolumn{2}{c|}{\textbf{Mean}} & \multicolumn{2}{c|}{\textbf{Variance}} & \multicolumn{2}{c|}{\textbf{Parameter} $a$} & \multicolumn{2}{c|}{\textbf{Parameter} $b$} & \multicolumn{1}{c|}{\textbf{Reference}} \\ 
		\hline 
		 $p_\mathrm{min}\sim\Gamma(a,b)$  & $5.7400$ & $[\mu\mathrm{g/g}]$ & $2.2464$ & $[\mu\mathrm{g^2/g^2}]$ & $13.667$ & $[-]$ & $2.5552$ & $[\mathrm{g}/\mu\mathrm{g}]$& {\cite{wu_decrease_2011} Fig 4C (AD)}  \\
		 $p_\Delta\sim\Gamma(a,b)$  & $3.0500$ & $[\mu\mathrm{g/g}]$ & $8.2143$ & $[\mu\mathrm{g^2/g^2}]$ & $0.1325$ & $[-]$ & $0.3713$ & $[\mathrm{g}/\mu\mathrm{g}]$& {\cite{wu_decrease_2011} Fig 4C (Control)}  \\
		 $q_\mathrm{max}\sim\Gamma(a,b)$  & $13.086$ & $[\mu\mathrm{g/g}]$ & $101.33$ & $[\mu\mathrm{g^2/g^2}]$ & $0.6884$ & $[-]$ & $0.1291$ & $[\mathrm{g}/\mu\mathrm{g}]$& {\cite{wu_decrease_2011} Tab 1 (AD)}  \\ \hline
	\end{tabular}
	\caption{Estimated distributions parameters for the amyloid-$\beta$ protein simulation \cite{wu_decrease_2011,roberts_biochemically-defined_2017}.}
	\label{tab:abeta_dist_param}
\end{table}
\section{Numerical methods}
\label{sec:methods}
In this section, after defining some preliminary concepts, we introduce all the methods adopted for discretizing problems~\eqref{eq:hm_strong} and~\eqref{eq:fk_strong}. We introduce a polytopic mesh partition $\partition$ of the domain $\Omega$ made of disjoint polygonal/polyhedral elements $K$. Moreover, we define the set of interfaces (namely $(d-1)-$dimensional facets of two neighboring elements) as $\faces$. This set is then decomposed into the union of interior faces ($\facesinternal$) and exterior faces ($\facesboundary$) lying on the boundary of the domain $\partial\Omega$, i.e. $\faces = \facesinternal \cup \facesboundary$. Concerning assumptions on the domain partition, we refer to the properties in \cite{dipietro:HHO,corti_discontinuous_2023}. 
\par
Let us define $\mathbb{P}_{\ell}(K)$ as the space of polynomials of total degree $\ell \geq 1$ over a mesh element $K$. Then we can introduce the following discontinuous finite element space $\Wh = \{w\in L^2(\Omega):\; w|_K\in\mathbb{P}_{\ell}(K)\quad\forall K\in\partition\}$. Moreover, we introduce the trace operators, namely the average operator $\averagel{\cdot}\averager$ and the jump operator $\jumpl{\cdot}\jumpr$ on $F\in\facesinternal$. For the complete definitions of the operators, we refer to \cite{arnoldUnifiedAnalysisDiscontinuous2001}. 
\subsection{PolyDG semi-discrete formulation of the heterodimer model}
To construct the semi-discrete formulation, we define the penalization functions
\begin{equation}
  \eta:\facesinternal\rightarrow\mathbb{R}_+ \mathrm{\;such\,that\;}\eta = \eta_0 
    \max\left\{\{d^K\}_\mathrm{H},\{k^K\}_\mathrm{H}\right\}\dfrac{\ell^2}{\{h\}_\mathrm{H}},
    \label{eq:hm_penalty_1}
\end{equation}
\begin{equation}
  \tilde{\eta}:\facesinternal\rightarrow\mathbb{R}_+ \mathrm{\;such\,that\;} \tilde{\eta} = \eta_0 
    \max\left\{\{d^K\}_\mathrm{H},\{\tilde{k}^K\}_\mathrm{H}\right\}\dfrac{\ell^2}{\{h\}_\mathrm{H}},
    \label{eq:hm_penalty_2}
\end{equation}
where $\eta_0$ is a constant parameter that should be chosen sufficiently large to ensure the stability of the discrete formulation, $d^K = \|\sqrt{\mathbf{D}|_K}\|^2$, $k^K = \|\, k_{12}|_K + k_1|_K \|$, and $k^K = \|\,k_{12}|_K + \tilde{k}_1|_K\|$. In Equations \eqref{eq:hm_penalty_1}and \eqref{eq:hm_penalty_2}, we are considering the harmonic average operator defined as $\{v\}_\mathrm{H} = \frac{2 v^+ v^-}{v^+ + v^-}$. We define the bilinear form $\mathcal{A}_h:\Wh\times \Wh\rightarrow \mathbb{R}$ as:
\begin{equation}
        \mathcal{A}_h(u,v) = \int_{\Omega} \left(\mathbf{D} \nabla_h u\right)\cdot\nabla_h v + \sum_{F\in\facesinternal}\int_{F}\left(\eta \jumpl u\jumpr \cdot \jumpl v\jumpr - \averagel\mathbf{D} \nabla_h u\averager \cdot \jumpl v \jumpr - \jumpl u\jumpr \cdot \averagel\mathbf{D} \nabla_h v \averager\right) \mathrm{d}\sigma \quad \forall u,v \in\Wh,
\label{eq:bilinearform}
\end{equation}
where $\nabla_h$ is the elementwise gradient. Analogously, we can define the bilinear form $\tilde{\mathcal{A}}_h$, by substituting in \ref{eq:bilinearform} the penalty $\eta$ with $\tilde{\eta}$. Given suitable discrete approximations $p_{0h},q_{0h}\in \Wh$ of the initial conditions of Equation~\eqref{eq:hm_strong}, the semi-discrete PolyDG formulation reads:
\par
\bigskip
For each $t>0$ and for each $\boldsymbol{v}\in\Gamma$, find $(p_h,q_h) =(p_h(t,\boldsymbol{v}),q_h(t,\boldsymbol{v}))\in \Wh\times\Wh$ such that:
\begin{equation}
{
\begin{cases}
  \left(\frac{\partial p_h}{\partial t}, v_h\right)_{\Omega} + \mathcal{A}_h(p_h,v_h) + (k_1\,p_h+k_{12}\,p_h q_h,v_h)_\Omega = (k_0,v_h)_\Omega & \forall v_h\in \Wh, \\
 \left(\frac{\partial q_h}{\partial t}, w_h\right)_{\Omega} + \tilde{\mathcal{A}}_h(q_h,w_h) + (\tilde{k}_1\,q_h-k_{12}\,p_h q_h,w_h)_\Omega = 0 & \forall w_h\in \Wh, \\
 p_h(\boldsymbol{x},0,\boldsymbol{v}) = p_{0h}(\boldsymbol{v}), \qquad q_h(\boldsymbol{x},0,\boldsymbol{v}) = q_{0h}(\boldsymbol{v}). 
\end{cases}}
\label{eq:hm_dgformulation}
\end{equation}
For details regarding the derivation and stability and \textit{a-priori} error estimates for this problem, we refer to \cite{antonietti_heterodimer_2023}.
\par
We can rewrite the obtained formulation in an algebraic form. We consider $\{\phi_j\}_{j=0}^{N_h}$ a set of basis function for $\Wh$, being $N_h$ its dimension. Then, we can expand the discrete solutions in terms of the chosen basis:
\begin{equation*}
    p_h(\boldsymbol{x},t,\boldsymbol{v}) = \sum_{j=0}^{N_h} \mathbf{P}_j(t,\boldsymbol{v}) \phi_j(\boldsymbol{x}) \quad
    q_h(\boldsymbol{x},t,\boldsymbol{v}) = \sum_{j=0}^{N_h} \mathbf{Q}_j(t,\boldsymbol{v}) \phi_j(\boldsymbol{x})
\end{equation*} 
We denote by $\mathbf{P}(t,\boldsymbol{v})$, $\mathbf{Q}(t,\boldsymbol{v})\in\mathbb{R}^{N_h}$, the corresponding vectors of the expansion coefficients. {Moreover, we define the following matrices $\mathrm{M},\mathrm{A},\mathrm{M}_\omega,\widehat{\mathrm{M}}_\omega\in\mathbb{R}^{N_h\times N_h}$ as:
\begin{align*}
     [\mathrm{M}]_{i,j}                 & = (\phi_j,\phi_i)_\Omega          & \mathrm{Mass\;matrix} \\
     [\mathrm{A}]_{i,j}                 & = \mathcal{A}_h(\phi_j,\phi_i)      & \mathrm{Stiffness\;matrix} \\
     [\mathrm{M}_\omega]_{i,j}               & = (\omega\phi_j,\phi_i)_\Omega              & \mathrm{Linear\;reaction\;matrix\;of\;a\;parameter\;}\omega \\
     [\widehat{\mathrm{M}}_\omega (\boldsymbol{\Theta})]_{i,j}               & = (\omega\Theta_h\phi_j,\phi_i)_\Omega              & \mathrm{Non-linear\;reaction\;matrix\;of\;a\;parameter\;}\omega\;\mathrm{and\,function}\;\Theta_h \\
\end{align*}
with $i,j = 1,...,N_h$, $\omega=\{\alpha,k_1,\tilde{k}_1,k_{12}\}$, $\Theta_h=\{c_h,q_h,p_h\}$ is the space-discrete multiplicative function in the nonlinear term and $\boldsymbol{\Theta}$ is the corresponding vector of the expansion coefficients. Finally, we introduce the right-hand side vector $[\mathbf{F}_{k_0}]_i= (k_0,\phi_i)_\Omega$  with $i = 1,...,N_h$.} Applying the definitions we can rewrite problem~\eqref{eq:hm_dgformulation} in algebraic form:
\par
\bigskip
For each $t>0$ and for each $\boldsymbol{v}\in\Gamma$, find $(\mathbf{P},\mathbf{Q}) = (\mathbf{P}(t,\boldsymbol{v}),\mathbf{Q}(t,\boldsymbol{v}))$ such that:
\begin{equation}
    \label{eq:HM_AlgebraicFormulation}
    \begin{cases}
    \mathrm{M}\dot{\mathbf{P}} + \mathrm{A}\mathbf{P} + \mathrm{M}_{k_1}\mathbf{P} + \widehat{\mathrm{M}}_{k_{12}}\left(\mathbf{P}\right) \mathbf{Q} =  \mathbf{F}_{k_0}, & \mathrm{in} \:  (0,T]\times\Gamma,\\[6pt] 
   \mathrm{M}\dot{\mathbf{Q}} + \mathrm{A}\mathbf{Q} + \mathrm{M}_{\tilde{k}_1}\mathbf{Q} - \widehat{\mathrm{M}}_{k_{12}}\left(\mathbf{Q}\right) \mathbf{P} =  \mathbf{0}, & \mathrm{in} \:  (0,T]\times\Gamma,\\[6pt]
   \mathbf{P}(0,\boldsymbol{v}) = \mathbf{P}_0(\boldsymbol{v}), \qquad \mathbf{Q}(0,\boldsymbol{v}) =\mathbf{Q}_0(\boldsymbol{v}), & \mathrm{in} \: \Gamma,
    \end{cases}
\end{equation}
where $\mathbf{P}_0$ and $\mathbf{Q}_0$ are the vector expansions associated to the initial conditions $p_{0h}$ and $q_{0h}$, respectively.
\subsection{PolyDG semi-discrete formulation of the Fisher-Kolmogorov model}
To construct the semi-discrete formulation, we define the penalization function
\begin{equation}
     \eta:\facesinternal\rightarrow\mathbb{R}_+ \mathrm{\;such\,that\;}\eta = \eta_0 \max\left\{\{d^K\}_\mathrm{H},\{\alpha^K\}_\mathrm{H}\right\}\dfrac{\ell^2}{\{h\}_\mathrm{H}},
    \label{eq:fk_penalty}
\end{equation}
where $\eta_0$ is a constant parameter that should be chosen sufficiently large to ensure the stability of the discrete formulation, $d^K$ is defined as before, and $\alpha^K = \|\alpha|_K\|$. 
\par
Using the definition of the penalty parameter in Equation~\eqref{eq:fk_penalty} in the definition of the bilinear form in Equation~\eqref{eq:bilinearform}, the semi-discrete PolyDG formulation reads:
\par
\bigskip
For each $t>0$ and for each $\boldsymbol{v}\in\Gamma$, find $c_h = c_h(t,\boldsymbol{v})\in \Wh$ such that:
\begin{equation}
\label{eq:fk_dgformulation}
    \begin{cases}
     \left(\dfrac{\partial c_h}{\partial t},w_h\right)_\Omega + \mathcal{A}_h(c_h,w_h) - (\alpha c_h,w_h)_\Omega + (\alpha c_h^2,w_h)_\Omega = 0
     \quad \forall w_h\in \Wh, \\
    c_h(\boldsymbol{x},0) = c_{0h}.
    \end{cases}
\end{equation}
where $c_{0h}\in\Wh$ is a suitable approximation of $c_0$. For details regarding the derivation of the semi-discrete formulation and the stability and \textit{a-priori} error analysis, we refer to \cite{corti_discontinuous_2023}. We consider the discrete solutions expansion in terms of the chosen basis:
\begin{equation*}
    c_h(\boldsymbol{x},t,\boldsymbol{v}) = \sum_{j=0}^{N_h} \mathbf{C}_j(t,\boldsymbol{v}) \phi_j(\boldsymbol{x}),
\end{equation*} 
and we denote by $\mathbf{C}(t,\boldsymbol{v})$ the corresponding vector of the expansion coefficients. Finally, the semi-discrete algebraic formulation of~\eqref{eq:fk_dgformulation} reads:
\par
\bigskip
For each $t>0$ and for each $\boldsymbol{v}\in\Gamma$, find $\mathbf{C} = \mathbf{C}(t,\boldsymbol{v})$ such that:
\begin{equation}
    \label{eq:FK_AlgebraicFormulation}
    \begin{cases}
   \mathrm{M}\dot{\mathbf{C}} + \mathrm{A}\mathbf{C} - \mathrm{M}_{\alpha}\mathbf{C} + \widehat{\mathrm{M}}_{\alpha}\left(\mathbf{C}\right) \mathbf{C} =  \mathbf{0}, & \mathrm{in} \:  (0,T]\times\Gamma,\\[6pt] 
   \mathbf{C}(0,\boldsymbol{v}) = \mathbf{C}_0(\boldsymbol{v}),  & \mathrm{in} \: \Gamma,\\[6pt]
    \end{cases}
\end{equation}
where $\mathbf{C}_0$ is the vector expansion associated to the initial conditions $c_{0h}$. 
\subsection{IMEX-PolyDG fully-discrete approximation of heterodimer model}
For the time discretization, we employ a linearly implicit IMEX Runge-Kutta method \cite{calvo_linearly_2001,hairer_solving_1996}. The advantages of these methods are the possibility of having a high-order approximation at the time of the problem and the explicit treatment of the nonlinear part of the equations, which avoids the necessity of nonlinear iterative solvers. {Let us introduce a generic differential problem continuous in time associated with the equation:
\begin{equation}
    \mathrm{M} \dot{\mathbf{y}}(t) = \mathbf{F}(t) + \mathrm{L}\mathbf{y}(t) + \mathrm{N}(\mathbf{y}(t)).
    \label{eq:genericIMEX}
\end{equation}
where $\mathbf{y}$ is the solution of the differential problem, $\mathrm{M}$ is the mass matrix obtained through the space discretization, $\mathrm{L}$ is the linear part of the differential problem, and $\mathrm{N}(\mathbf{y})$ is the nonlinear one.} We introduce the following Butcher tableaux for the time discretization of the linear and nonlinear parts, respectively:
\begin{equation}
\begin{array}
{c|ccccc}
c_1 & \gamma \\
c_2 & a_{21} & \gamma \\
\vdots & \vdots & \ddots & \ddots \\
c_{s} & a_{s\,1} & \dots & a_{s\,s-1} & \gamma \\
\hline
& b_1 & \dots & b_{s-1} & \gamma 
\end{array}
\qquad
\begin{array}
{c|ccccc}
& 0 \\
& \hat{a}_{21} & 0 \\
& \hat{a}_{31} & \hat{a}_{32} & 0 \\
& \vdots & \vdots & \ddots & \ddots \\
\gamma & \hat{a}_{s+1\,1} & \hat{a}_{s+1\,2} & \dots & \hat{a}_{s+1\,s} & 0 \\
\hline
 \gamma & \hat{b}_1 & \hat{b}_2 & \dots & \hat{b}_s & \hat{b}_{s+1} 
\end{array}.
\end{equation}
We construct a partition $0 < t_1 < t_2 < ... < t_{N_T} = T$ of the time interval $[0,T]$ into $N_T$ intervals of constant time step $\Delta t = t_{n+1}-t_n$. Then the time discretization of problem \eqref{eq:genericIMEX} can be written as:
\par
\bigskip
For $n=0,...,N_T-1$ and given $\mathbf{y}^n$, find $\mathbf{y}^{n+1}$ such that:
\begin{equation}
    \begin{dcases}
        \mathrm{M} \mathbf{K}_1 = N(\mathbf{y}^n), \\
        (\mathrm{M}-\Delta t \gamma\mathrm{L})\mathbf{K}_i = \mathbf{F}(t^n+c_j \Delta t) + L\left(\mathbf{y}^n+\Delta t \left(\sum_{j=1}^{i-1} a_{ij} \mathbf{K}_j+\sum_{j=1}^i \hat{a}_{ij} \hat{\mathbf{K}}_j\right)\right), & i=1,...,s-1, \\
        \bar{\mathbf{y}}^i = \mathbf{y}^n + \Delta t \sum_{j=1}^{i}\left( a_{ij} \mathbf{K}_j+ \hat{a}_{ij} \hat{\mathbf{K}}_j\right), & i=1,...,s-1, \\
        \mathrm{M} \mathbf{K}_{i+1} = N(\bar{\mathbf{y}}^i), & i=1,...,s-1,\\
        \mathbf{y}^{n+1} = \mathbf{y}^n + \Delta t \left(\sum_{j=1}^{s} b_{j} \mathbf{K}_j  + \sum_{j=1}^{s+1} \hat{b}_{j} \hat{\mathbf{K}}_j\right). \\
    \end{dcases}
\end{equation}
Applying this multistage method to the two problems allows us to obtain an IMEX Runge-Kutta time-stepping method with a level of accuracy that can be arbitrarily high with appropriate choices of coefficients. We refer to \cite{calvo_linearly_2001,hairer_solving_1996} for examples of tableaux with different orders of convergence.
\subsection{Mesh construction}
\begin{figure}[t]
     \centering
     \begin{subfigure}[b]{0.23\textwidth}
         \centering
         \includegraphics[width=\textwidth]{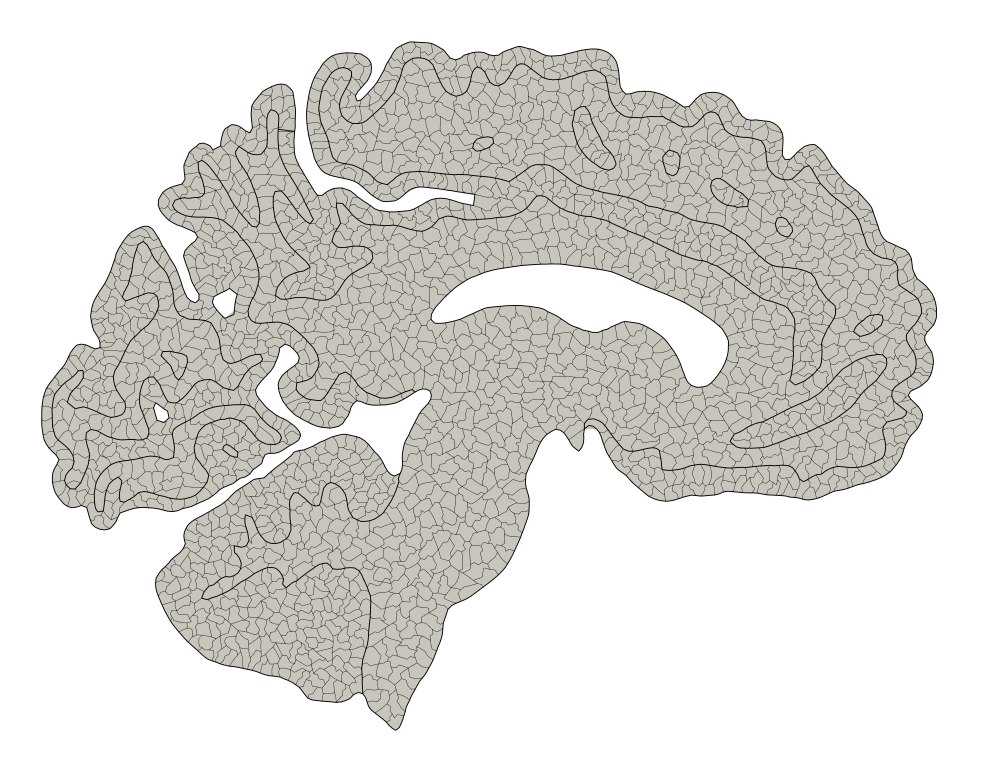}
         \caption{Polygonal mesh}
         \label{fig:BrainSection:Mesh}
     \end{subfigure}
     \hfill
     \begin{subfigure}[b]{0.235\textwidth}
         \centering
         \includegraphics[width=\textwidth]{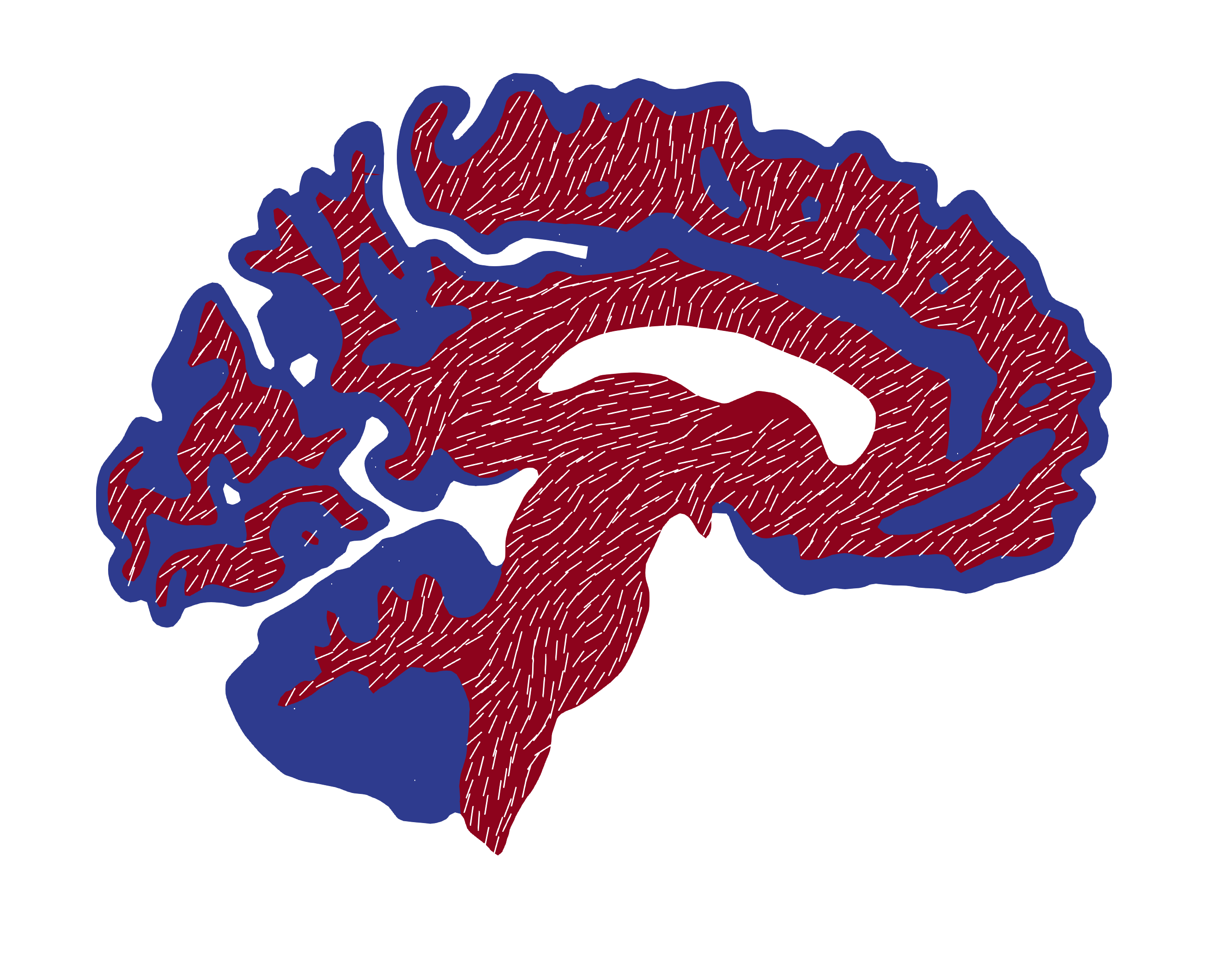}
         \caption{Axonal directions}
         \label{fig:BrainSection:Axn}
    \end{subfigure}
    \hfill
    \begin{subfigure}[b]{0.255\textwidth}
         \centering
         \includegraphics[width=\textwidth]{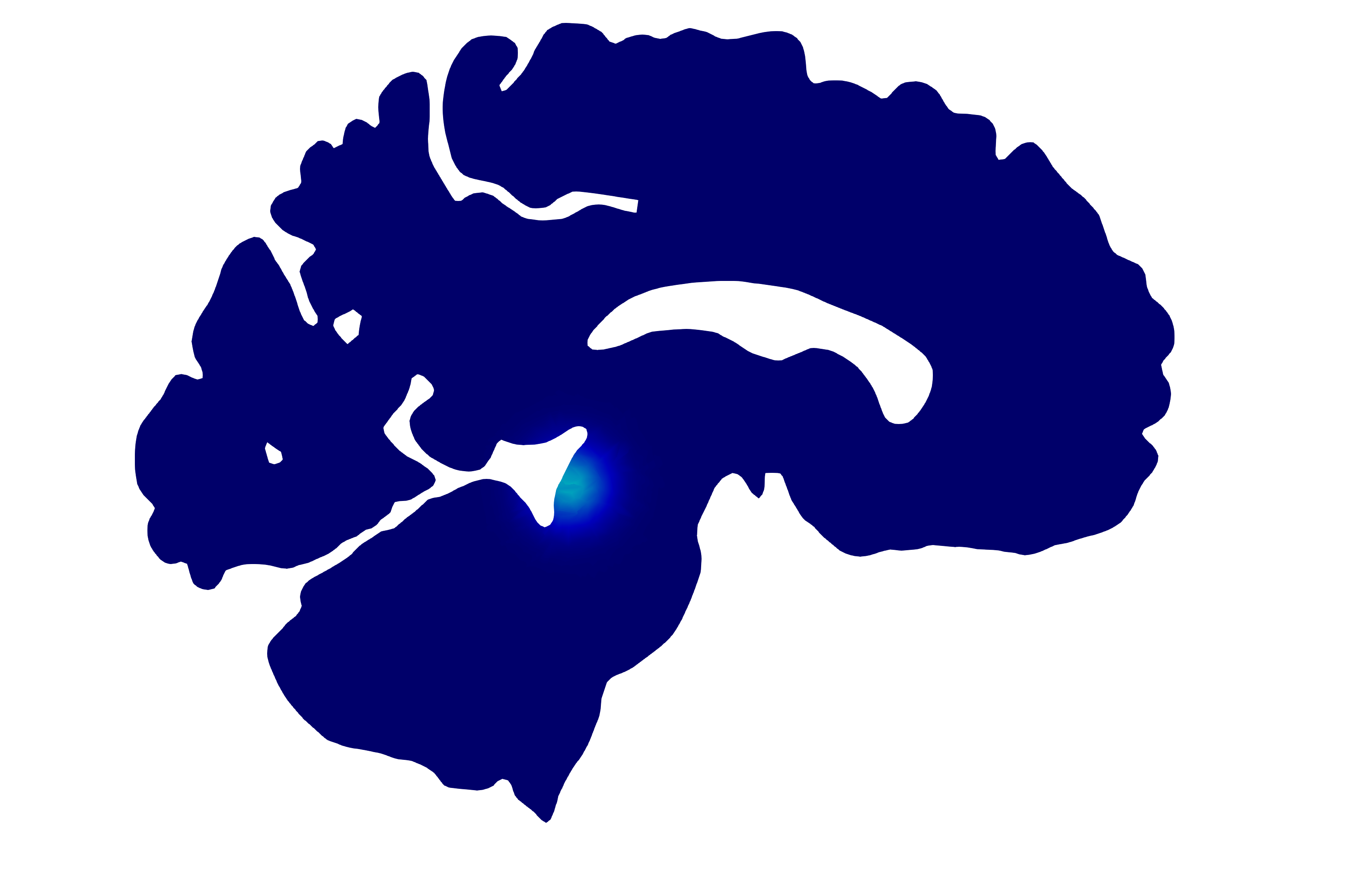}
         \caption{Tau protein $(t=0)$}
         \label{fig:BrainSection:tau}
    \end{subfigure}
    \hfill
    \begin{subfigure}[b]{0.255\textwidth}
         \centering
         \includegraphics[width=\textwidth]{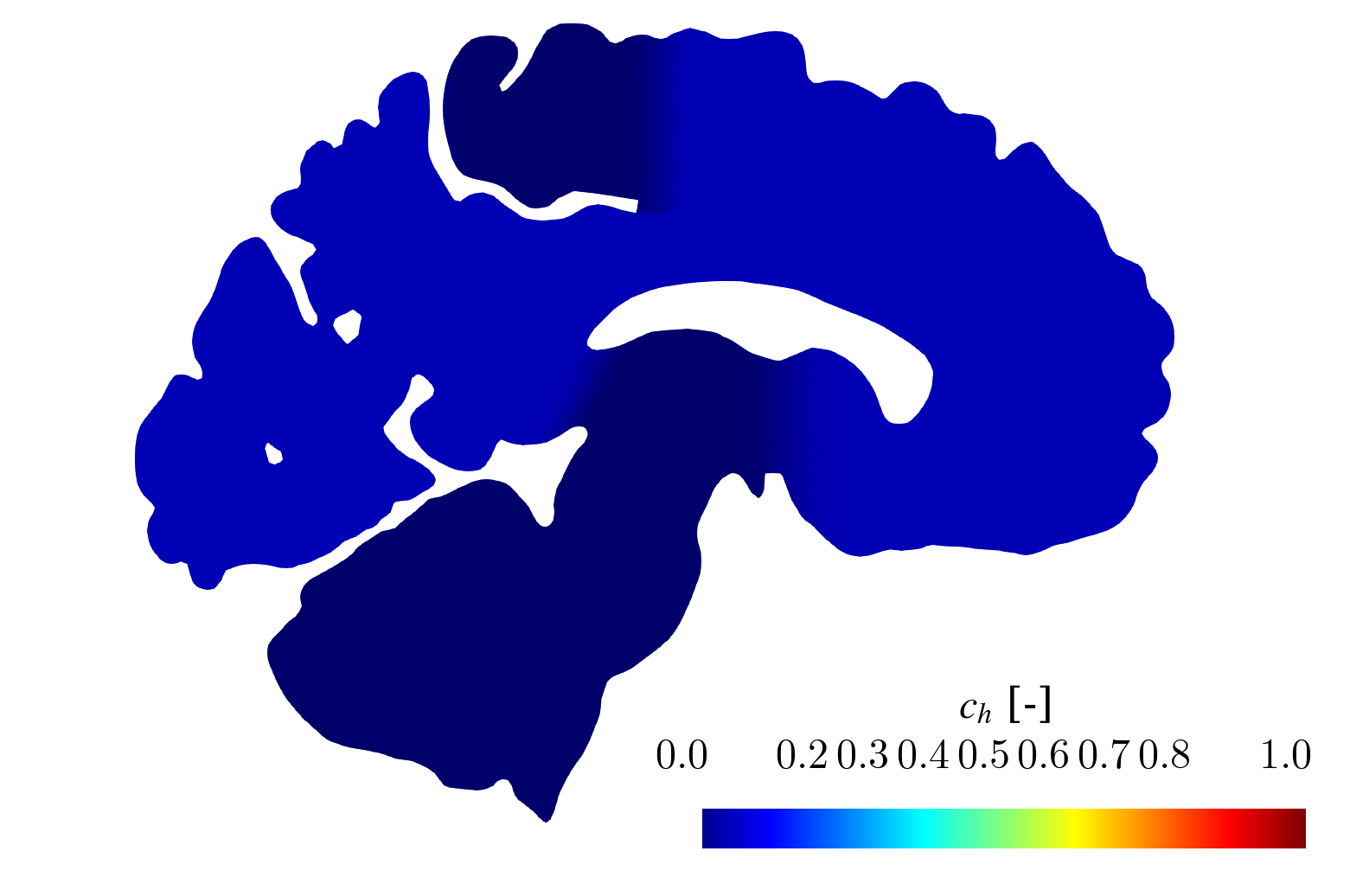}
         \caption{Amyloid-$\beta$ $(t=0)$}
         \label{fig:BrainSection:beta}
    \end{subfigure}
        \caption{Brain section with the specification of the polygonal mesh (a), the axonal directions with the distinction between white (red) and grey (blue) matters (b), and the initial conditions of the simulation of tau protein (c) and amyloid-$\beta$, respectively.}
        \label{fig:BrainSection}
\end{figure}

The numerical simulations are developed using the open source \texttt{lymph} library \cite{antonietti_lymph_2024}, implementing the PolyDG method for multiphysics problems. For the simulations, we consider a mesh of a sagittal brain section obtained by segmenting a structural magnetic resonance image (MRI) from the OASIS-3 database \cite{OASIS3}. We perform the segmentation using Freesurfer \cite{Freesurfer}. We construct an initial detailed triangular grid of $43\,402$ triangles. Then, we agglomerate the latter one employing ParMETIS \cite{Parmetis}. The final mesh comprises $2\,031$ polygonal elements, as shown in Figure~\ref{fig:BrainSection:Mesh}. 
\par
Constructing the white matter's axonal directions is fundamental for defining the diffusion tensor $\mathbf{D}$. These can be derived from Diffusion-Weighted Images (DWI) and computing the principal eigenvector of the diffusion tensor in each voxel of the image through Nibabel \cite{Nibabel}. Concerning the complete procedure, we refer to \cite{Mardal:Mesh}. Figure~\ref{fig:BrainSection:Axn} reports the resulting axonal direction.
\par
Concerning the physical parameters, which are not the object of the sensitivity analysis, we fix the conversion rate $k_{12} = 0.2 \mathrm{years}^{-1}$ as in \cite{antonietti_heterodimer_2023}. The extracellular diffusion is imposed to be $d_\mathrm{ext}=8\times10^{-6}$, while the axonal diffusion is set $d_\mathrm{ext}=8\times10^{-5}$ in the white matter, and it is $d_\mathrm{ext}=0$ in the gray matter \cite{corti_discontinuous_2023}. For what concerns the discretization parameters in space, we fix the polynomial order $\ell=5$ in every element of the discretization, and we set the penalty parameter $\gamma_0=10$. Finally, we adopt a time step $\Delta t = 0.025$ and simulate up to $T=40$ years to have a complete sensitivity analysis over a time interval larger than the typical scales of Alzheimer's disease \cite{van_oostveen_imaging_2021}.
\par
The initial conditions are computed according to other literature works associated with the same proteins \cite{weickenmeier_multiphysics_2018}. In particular, for the tau protein-spreading simulation, we initially locate the misfolding proteins in the entorhinal cortex \cite{braak_alzheimer_1991}, as in Figure~\ref{fig:BrainSection:tau}. For what concerns the amyloid-$\beta$ protein, we impose initial misfolding protein concentrations in the neocortex \cite{jucker_pathogenic_2011}, as reported in Figure~\ref{fig:BrainSection:beta}.
\section{Results}
\label{sec:results}
In this section, we report the sensitivity analysis results of the models applied to spreading the protein tau and the amyloid-$\beta$ in the brain. In all the analysis, we will also present the space average of the solutions, computed for a generic function $g(\boldsymbol{x})$ as:
\begin{equation}
    \langle g(\boldsymbol{x}) \rangle = \frac{1}{|\Omega|}\int_\Omega g(\boldsymbol{x}) \mathrm{d}\boldsymbol{x}.
\end{equation}
\subsection{Sensitivity analysis in tau protein spreading}
In this section, we present the results associated with the application of tau protein spreading in Alzheimer's disease, with the parameters calibrated starting from the distributions presented in Table \ref{tab:tau_dist_param} and derived from biological measurements in \cite{wu_decrease_2011}.
\subsubsection{Sensitivity of heterodimer model} 
\begin{figure}[t!]
    \begin{subfigure}[b]{\textwidth}
         \centering
        \includegraphics[width=0.95\textwidth]{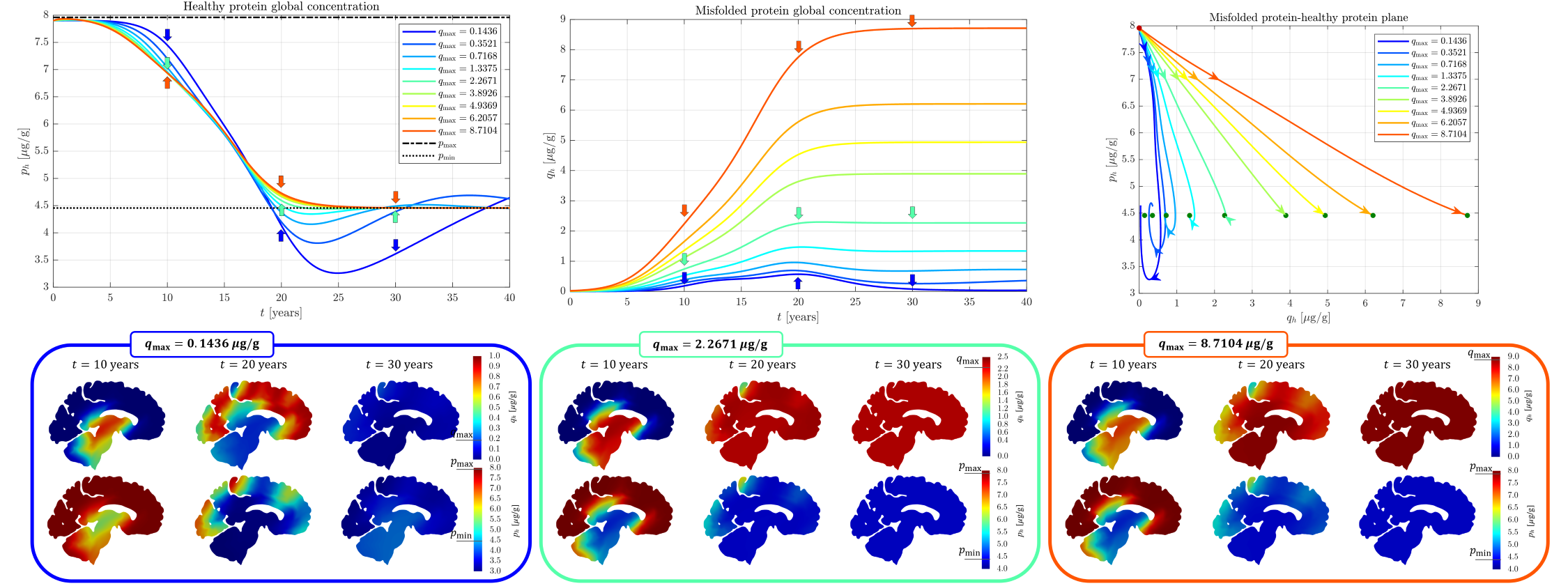}
        \caption{Sensitivity to $q_\mathrm{max}$ of heterodimer in tau protein spreading}
        \label{fig:sens_tau_qmax}
    \end{subfigure}
    \begin{subfigure}[b]{\textwidth}
         \centering
        \includegraphics[width=0.95\textwidth]{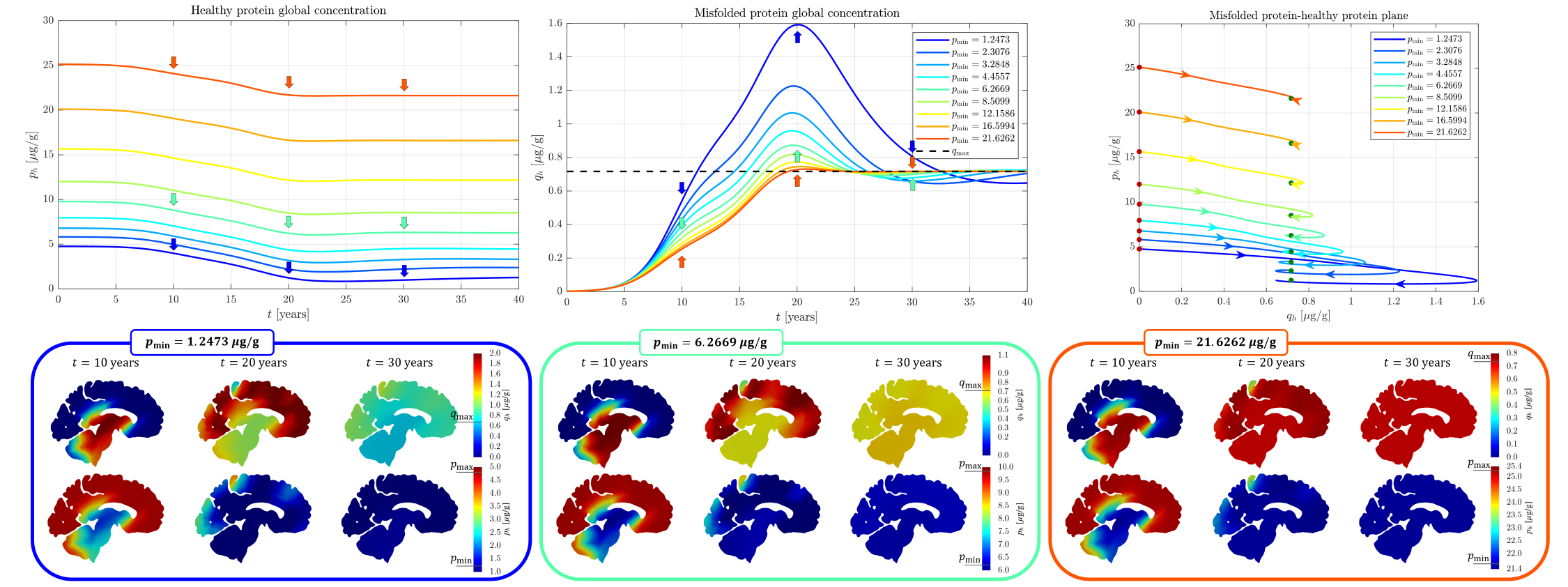}
        \caption{Sensitivity to $p_\mathrm{min}$ of heterodimer in tau protein spreading}
        \label{fig:sens_tau_pmin}
    \end{subfigure}
    \begin{subfigure}[b]{\textwidth}
         \centering
        \includegraphics[width=0.95\textwidth]{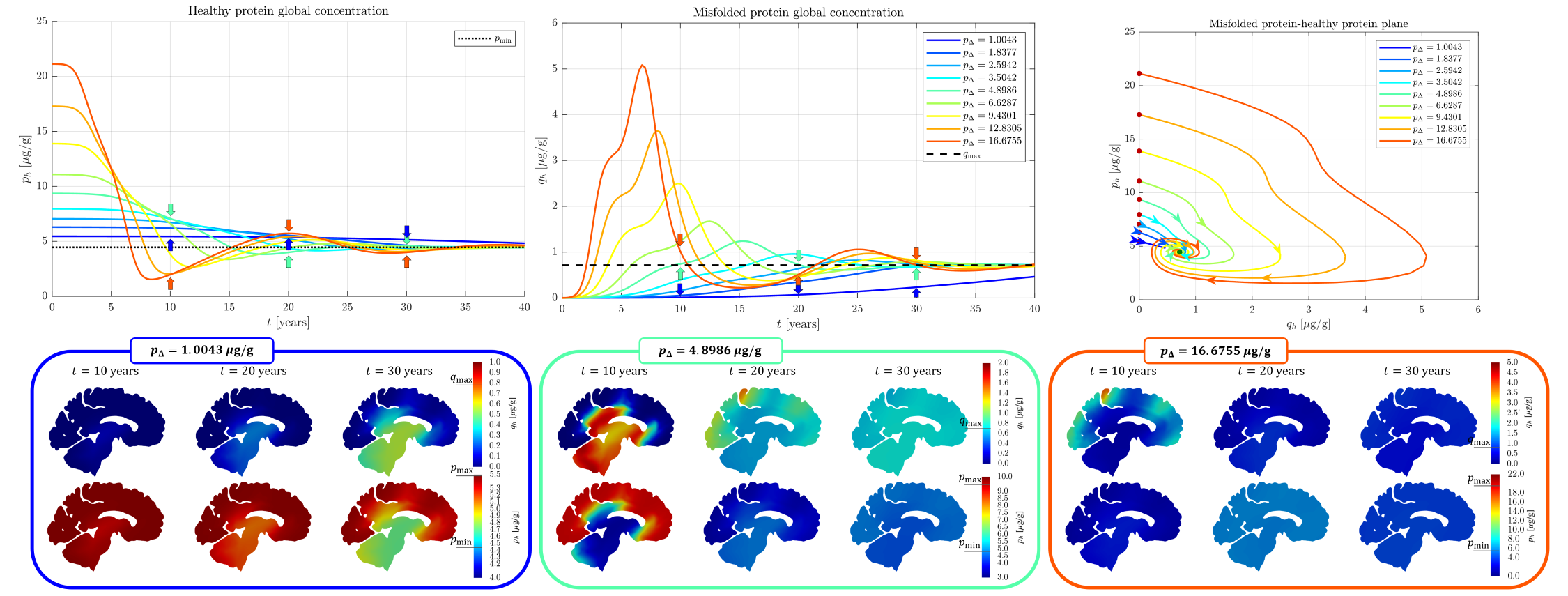}
        \caption{Sensitivity to $p_\Delta$ of heterodimer in tau protein spreading}
        \label{fig:sens_tau_pdelta}
    \end{subfigure}     
   \caption{Sensitivity analysis of heterodimer model in tau protein for the variation of the parameter $q_\mathrm{max}$ (a), $p_\mathrm{min}$ (b), $p_\Delta$ (c). For each parameter, in the first row, we report the space average of the healthy protein in time (left), of the misfolded protein in time (center), and of both in the phase space (right). In the second row, we report some snapshots of the solutions for three selected parameter values.}
   \label{fig:sens_tau}
\end{figure}
The analysis shows a change in the behavior of the solution unless for each value of the parameter $q_\mathrm{max}$, the stability of the equilibria does not change. However, it is visible from Figure~\ref{fig:sens_tau_qmax} that the equilibrium is a stable focus for small parameter values. Indeed, the dynamics of the solutions in the proximity of the value present oscillations. The equilibrium undergoes a bifurcation that transforms it into a stable node, increasing the value \cite{arnold2013dynamical}. {Coherently with Figure \ref{fig:bifsurf}, the bifurcation value can be estimated to be $q_\mathrm{max}\simeq3.4541\,\mu\mathrm{g/g}$.} Considering the probability distribution of the parameter $q_\mathrm{max}$, the oscillatory behavior is associated with the region of most probable values, and it is also present in the simulation with the average value $q_\mathrm{max}=0.7168 \,\mu\mathrm{g/g}$. The oscillations are also visible spatially, watching at the lower line of Figure~\ref{fig:sens_tau_qmax}.
\par
The sensitivity analysis of the solutions for variations of the equilibrium concentration of healthy proteins $p_\mathrm{min}$ shows that equilibrium is a stable focus for each parameter value. It is visible from Figure~\ref{fig:sens_tau_pmin}, where it can be observed that we have more significant oscillations in the dynamics of the space average of misfolded proteins $q$ for small parameter values. Indeed, a decrease in the parameter value causes a sensible reduction of the real part of the eigenvalues, making the oscillating component (associated with the imaginary part) dominant. Looking at the lower line of Figure~\ref{fig:sens_tau_pmin}, we can observe the impact of the oscillation on the solution fields only for the small values of $p_\mathrm{min}$, {suggesting that for large values of total tau in a patient, the longitudinal oscillations on this biomarker can be not detectable.}
\par
The analysis of the impact of the variation of healthy proteins $p_\Delta$ on the problem solutions shows a stable focus behavior for most of the tested values. For high values of the parameter $p_\Delta$, the oscillations become more prominent, with an initial fast increase of the population of misfolded proteins, due to the high presence of the healthy ones, as visible in Figure~\ref{fig:sens_tau_pdelta}. {This fact would suggest that patients who undergo fast hyperphosphorylation of the tau-protein can be associated with oscillating behavior in the concentrations if tested longitudinally. It can be observed that for small values of the parameter, we do not notice the oscillations, indeed the equilibrium is a stable node for values of $p_\Delta$ smaller than $1.1177\,\mu\mathrm{g}/\mathrm{g}$ (see Figure~\ref{fig:bifsurf}).} Finally, in the lower line of Figure~\ref{fig:sens_tau_pdelta}, we can observe the impact of the oscillation on the solution fields and confirm the absence of oscillations for the smallest one.
\subsubsection{Sensitivity of FK model} 
\begin{figure}[t]
    \centering
    \includegraphics[width=\textwidth]{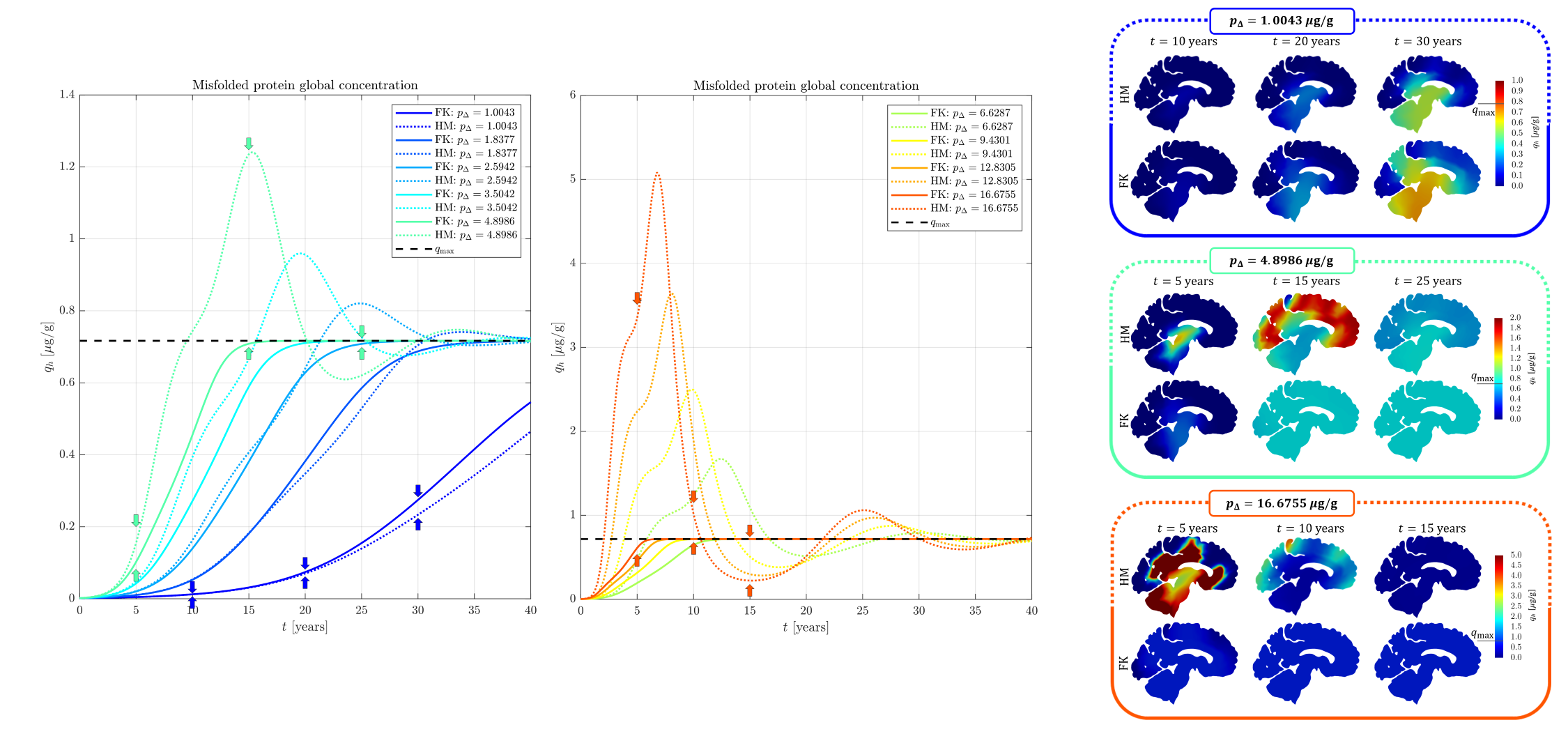}
   \caption{Sensitivity analysis of the FK model in tau protein for the variation of the parameter {$ p_\Delta$}. We compare the analysis results between FK and heterodimer in two different graphs. In the third column, we report some snapshots of the solutions for three selected parameter values.}
   \label{fig:sens_tau_FK}
\end{figure}
The sensitivity analysis concerning the value of $p_\Delta$ can also be performed on the FK model. Indeed, the reaction parameter $\alpha$ depends only on $p_\Delta$ and not on $p_\mathrm{min}$ and $q_\mathrm{max}$. The analysis shows that by increasing the initial concentration of healthy proteins, the development of the misfolded protein is faster. Indeed, looking at Figure~\ref{fig:sens_tau_FK}, we can notice this effect. To compare with the heterodimer model, we report in Figure~\ref{fig:sens_tau_FK} a non-normalized value $q_h=q_\mathrm{max}\,c_h$ for the FK solution. We can observe that because the FK model is composed of a single equation, the oscillatory behavior obtained with the heterodimer model for most of the parameter values is lost \cite{smoller_shock_1994}. Figure~\ref{fig:sens_tau_FK} shows large differences in the space-averaged solutions for large values of the parameter, while for the small ones, the solutions of the two models are comparable. For small values of the parameter, in which the stable equilibrium of the heterodimer is a node, the FK model generates a faster pathology.
\subsection{Sensitivity analysis in amyloid beta spreading}
In this section, we present the results associated with the application of amyloid-$\beta$ spreading in Alzheimer's disease, with the parameters calibrated starting from the distributions presented in Table \ref{tab:abeta_dist_param} and derived from biological measurements in \cite{wu_decrease_2011,roberts_biochemically-defined_2017}.
\subsubsection{Sensitivity of heterodimer model}
\begin{figure}[t!]
    \begin{subfigure}[b]{\textwidth}
         \centering
        \includegraphics[width=0.97\textwidth]{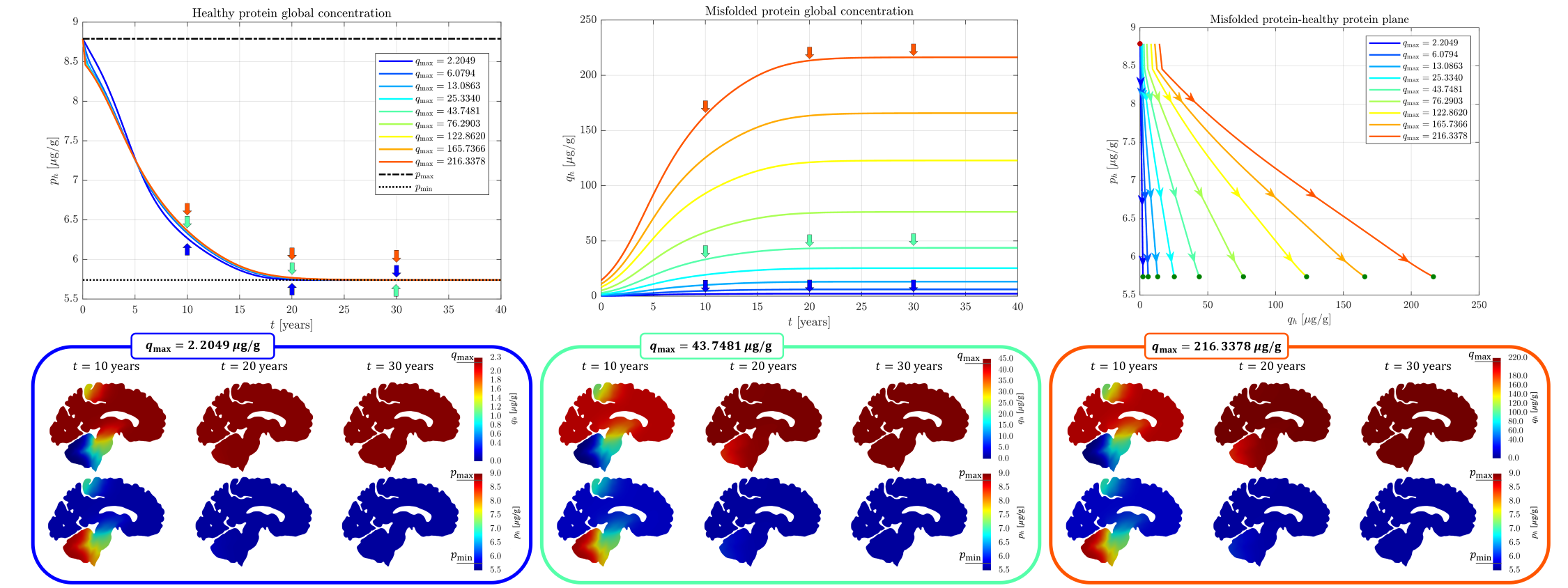}
        \caption{Sensitivity to $q_\mathrm{max}$ of heterodimer in amyloid-$\beta$ spreading}
        \label{fig:sens_abeta_qmax}
    \end{subfigure}
    \begin{subfigure}[b]{\textwidth}
         \centering
        \includegraphics[width=0.97\textwidth]{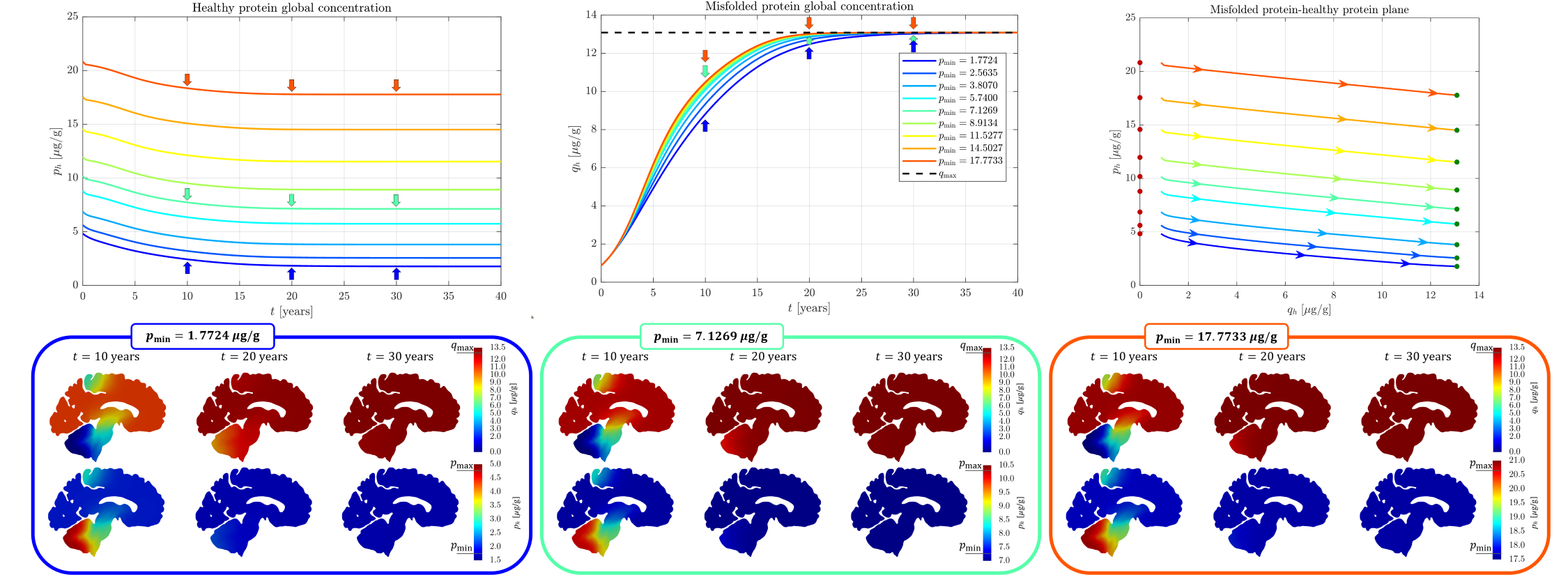}
        \caption{Sensitivity to $p_\mathrm{min}$ of heterodimer in amyloid-$\beta$ spreading}
        \label{fig:sens_abeta_pmin}
    \end{subfigure}
    \begin{subfigure}[b]{\textwidth}
         \centering
        \includegraphics[width=0.97\textwidth]{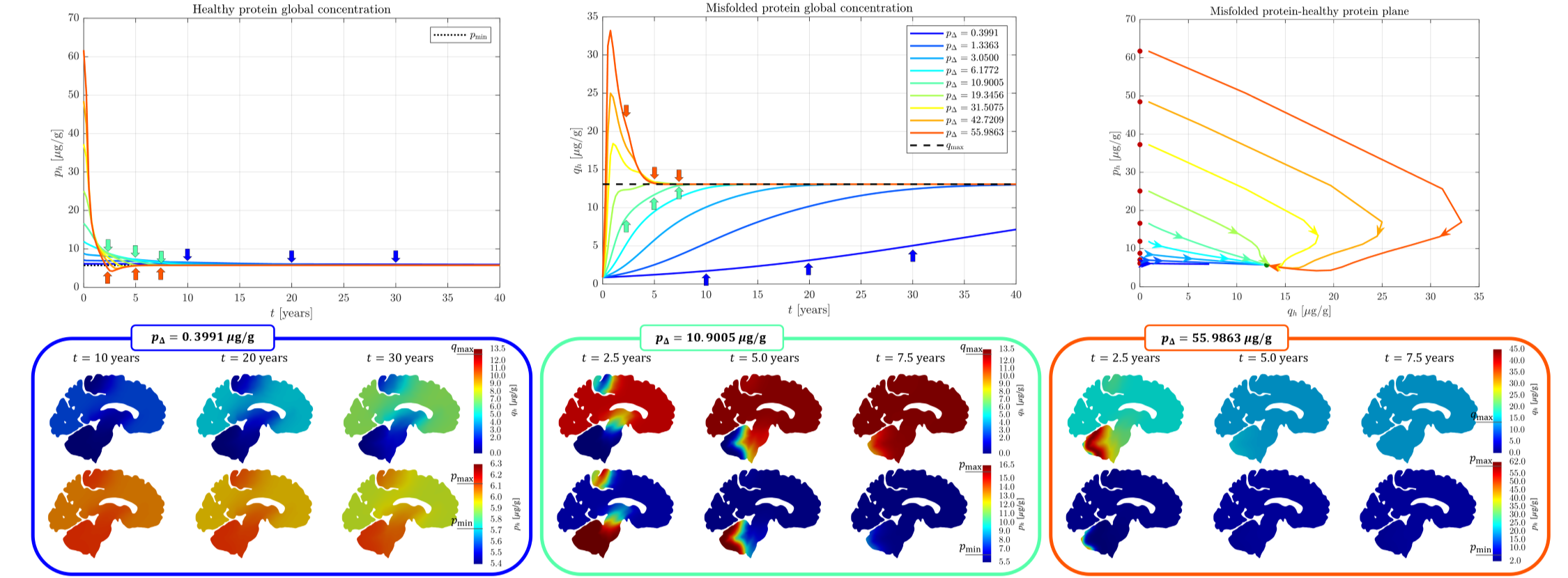}
        \caption{Sensitivity to $p_\Delta$ of heterodimer in amyloid-$\beta$ protein spreading}
        \label{fig:sens_abeta_pdelta}
    \end{subfigure}     
   \caption{Sensitivity analysis of heterodimer model in amyloid-$\beta$ for the variation of the parameter $q_\mathrm{max}$ (a), $p_\mathrm{min}$ (b), $p_\Delta$ (c). For each parameter, in the first row, we report the space average of the healthy protein in time (left), of the misfolded protein in time (center), and of both in the phase space (right). In the second row, we report some snapshots of the solutions for three selected parameter values.}
   \label{fig:sens_abeta}
\end{figure}
The analysis shows that the parameter $q_\mathrm{max}$ predominantly influences misfolded proteins' dynamics, while the impact on the healthy protein dynamics is almost inappreciable. Concerning the stable equilibrium, unless not visible in the range of the values computed for the amyloid-$\beta$, a bifurcation occurs at the value $q_\mathrm{max}=2.7644\,\mu\mathrm{g}/\mathrm{g}$ {(see Figure \ref{fig:bifsurf})}. Indeed, for values smaller than this, the equilibrium is a stable focus, and then it becomes a stable node. This fact is not visible in Figure~\ref{fig:sens_abeta_qmax}, where the oscillations are minimal for the value under this threshold. Indeed, the change of magnitude of the protein does not particularly affect the qualitative development of the pathology.
\par
The sensitivity analysis of the solutions to the variations of the parameter $p_\mathrm{min}$ shows that in the considered range of values, the stable equilibrium is a mode{, as visible in Figure \ref{fig:bifsurf}.} The impact of increasing the concentrations of healthy protein is a slight acceleration in the disease progression, which reaches a high level of misfolded proteins ($\simeq 90\%$ of $q_\mathrm{max}$) in 13-17 years. Looking at the lower line of Figure~\ref{fig:sens_abeta_pmin}, we can notice the acceleration at the time $t=10$ years. 
\par
The analysis of the impact of the variation of healthy proteins $p_\Delta$ on the problem solutions shows absent or small oscillations for most of the analyzed values. These are present only for high values of the parameter $p_\Delta$, with an initial fast increase of the population of misfolded proteins due to the high presence of the healthy ones, as visible in Figure~\ref{fig:sens_abeta_pdelta}. {The estimated bifurcation value from node to focus is $p_\Delta = 17.6845\,\mu\mathrm{g}/\mathrm{g}$.} Finally, in the lower line of Figure~\ref{fig:sens_abeta_pdelta}, we can observe the impact of the oscillation on the solution fields only for the largest of the three showed values.
\par
{Finally, in most of the cases and in particular for the most probable values of the parameters according to the reconstructed distributions, we do not expect oscillations in the concentrations of amyloid-$\beta$ levels in the cortex.}
\subsubsection{Sensitivity of FK model}
\begin{figure}[t]
    \centering
    \includegraphics[width=\textwidth]{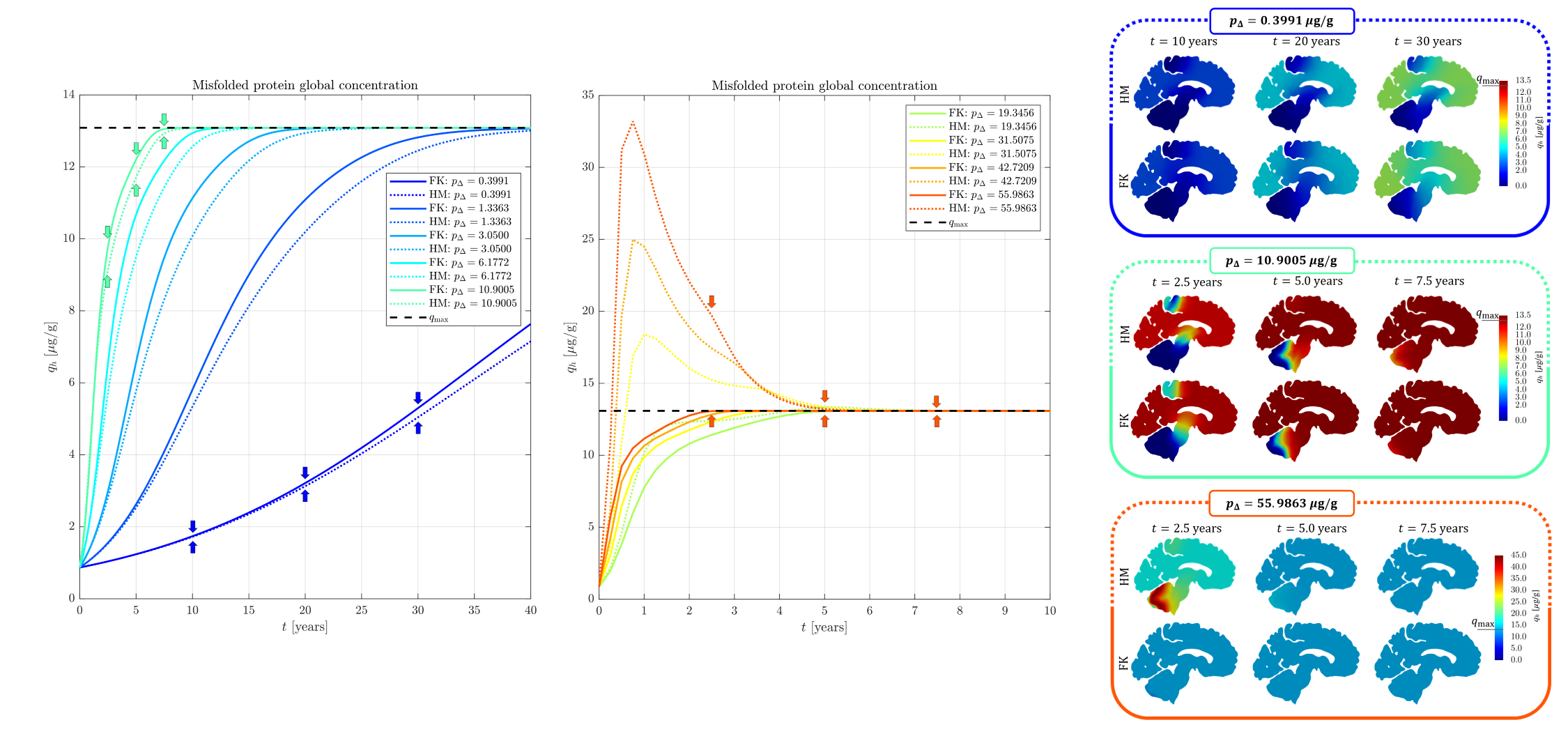}
   \caption{Sensitivity analysis of the FK model in amyloid-$\beta$ for the variation of the parameter {$ p_\Delta$}. We compare the analysis results between FK and heterodimer in two different graphs. In the third column, we report some snapshots of the solutions for three selected parameter values.}
   \label{fig:sens_beta_FK}
\end{figure}
We perform the latter parameter's sensitivity analysis on the FK model. As we have done in the tau protein case, we analyze the values of $q_h=q_\mathrm{max}\,c_h$ for the FK solution. The analysis shows that the concentration of misfolded protein increases faster with an increase of $p_\Delta$ and consequently of the parameter $\alpha$. In Figure~\ref{fig:sens_beta_FK}, we also report a comparison with the heterodimer model solution, and we can observe that the two models generate comparable solutions for small values of $p_\Delta$. We obtain different behaviors for the largest values of the parameter due to the focused nature of the stable equilibrium in the heterodimer model, causing sensible oscillations. This is also visible in the last column of Figure~\ref{fig:sens_beta_FK} for $p_\Delta = 55.9863\,\mu\mathrm{g}/\mathrm{g}$, which shows large oscillations in the solutions at the first time instant ($2.5$ years) in the heterodimer model, which are not obtainable with the FK model.
\section{Discussion}
\label{sec:discussion}
This section discusses the results obtained from the sensitivity analysis performed in section \ref{sec:results}. In particular, we discuss in section \ref{sec:disc_tau} the results associated with tau protein and in section \ref{sec:disc_abeta} the ones of amyloid-$\beta$. Finally, in section \ref{sec:limitations}, we report the study's limitations.
\subsection{Sensitivity analysis in tau protein spreading}
\label{sec:disc_tau}
The sensitivity analysis of the heterodimer model applied to the spreading of tau protein shows that the stable equilibrium of the system is a focus, for most of the tested values. This nature causes an oscillatory behavior of the space-averaged solutions and the solution fields locally, as reported in Figure \ref{fig:sens_tau}. Those oscillations are visible for values of the parameters that are derived from the medical literature (see Figure 5 in \cite{wu_decrease_2011}). {However, we underline that the largest oscillations are connected with initial values of total tau concentrations that are uncommon because associated with the tail of the Gamma distribution.}
\par
The oscillatory nature of these concentrations is not described in the literature. {Since these quantities are directly measurable only with histological \textit{post-mortem} examinations \cite{wu_decrease_2011}, we do not have longitudinal data to compare with our simulation. Literature findings report in a few cases fluctuating behavior of both the concentrations of tau protein in phosphorylated and healthy forms (see Figure 1H-G \cite{mcdade_longitudinal_2018}). However, the variation in the CSF measurement cannot be directly associated with the variation in the global cortical measure of the tau protein. Indeed, this can be related to a failure of the clearance mechanisms, due to the formation of large agglomerates of phosphorylated tau. Moreover, fluctuations in the CSF concentrations of tau \cite{moghekar_cerebrospinal_2012} are reported at the time scale of the circadian cycle, which complicates the evaluation of a connection between cortical and CSF concentrations. Concerning the cortical PET measurement, some localized decrease in SUVR measures seems possible (see Figure 3D-F \cite{sanchez_longitudinal_2021}). This can be coherent with the spatial oscillation detected by our model, however, the yearly oscillations are small and uncommon and can be associated with variations in the scanner parameters.}
\par
{Another interpretation concerning the absence of significant medical literature about oscillatory behavior can be associated with the time scales. Indeed, considering the mean values of the distribution, the beginning of the decreasing phase is around 15-20 years from the tau accumulation beginning, and these values are near to the total life expectancy of the diseased patients \cite{scheltens_alzheimers_2021}. An interesting result of the simulations in Figure \ref{fig:sens_tau}, is a two-stage behavior in the accumulation of tau protein, that slows down after a fast increase and then rapidly increases in a second time. This behavior has been recently found in the analysis of tau accumulation in Alzheimer's disease \cite{stockman_two_2024}.}
\par
{Sensible variations in the healthy tau concentrations do not allow a correct simplification of the heterodimer to the FK model \cite{weickenmeier_multiphysics_2018}. Indeed, the comparative analysis in Figure \ref{fig:sens_tau_FK} reports that the FK model cannot describe the fluctuations detected by the heterodimer in tau protein spreading. These findings suggest that the choice of the correct model should be fundamental in the simulation of tauopathies to accurately predict the time advancement of the pathology. However, it is unclear if the presence of oscillation is physically meaningful and which is the best model in this description.}
\par
{Looking at the medical interpretation of the sensitivity analysis of both models, it can be observed in Figure \ref{fig:sens_tau_FK} that individuals with high levels of healthy tau protein at the beginning of the pathology are associated with a faster misfolding process and with a faster decline of the pathology. In particular, this effect is amplified in the heterodimer simulation due to the focus behavior of the equilibrium. Literature findings report a correlation between high tau levels in CSF and a faster decline in dementia, supporting the idea that patients with naturally high tau levels can undergo severe dementia \cite{degerman_gunnarsson_high_2014}.}

\subsection{Sensitivity analysis in amyloid-$\beta$ spreading}
\label{sec:disc_abeta}
Concerning the application of the heterodimer model to the amyloid-$\beta$ spreading, the parameters are calibrated according to Table \ref{tab:abeta_dist_param}. In Figure \ref{fig:sens_abeta}, we can notice that we cannot detect any oscillation in either the misfolded or the healthy protein concentrations for most of the analyzed parameters. {We observe bifurcations only for extremely high values of $p_\Delta$, associated with the tail of the Gamma distribution, and so less probable and far from the average values computed from analysis in \cite{wu_decrease_2011,roberts_biochemically-defined_2017}.}
\par
{The results obtained in the sensitivity analysis suggest that the sigmoidal shape typically presented in the literature of the amyloid-$\beta$ accumulation is well caught by the model in most cases \cite{jagust_temporal_2021}. The impact of the parameters $q_\mathrm{max}$ and $p_\mathrm{min}$ in the velocity of the disease development seems to be not fundamental. Indeed, the increase of the $q_h/q_\mathrm{max}$ curve (Figure \ref{fig:sens_abeta}) is slightly affected by the changes in the parameters. The parameter $p_\Delta$ (namely the percentage of APP loss) is the most important factor to analyze the development of the pathology. In particular, the large values produce oscillations that indicate an initial accumulation of amyloid-$\beta$ that did not have enough time to be compensated by clearance in the first stages of the pathology. However, we underline that this effect is not reported in the medical literature.}
\par
{Some oscillations in the Florbetapir PET measurement have been reported in literature \cite{jagust_temporal_2021}. However, the study did not correlate these results with temporary decreases in the amyloid-$\beta$ accumulation, which is confirmed to have a sigmoidal longitudinal increment. Also \cite{mcdade_longitudinal_2018} (Figure 1E) shows only in a few cases decreased values in the cortical amyloid-$\beta$ PET. The article reports some oscillations in the cognition tests, but they seem to be not correlated with decreases in the amyloid-$\beta$ levels. Moreover, we remember that literature reports fluctuations in the attention tests in Alzheimer's diseased patients \cite{ballard_attention_2001} at the time scale of the circadian cycle, which needs to be taken into account in the complete picture of connecting cognition to biomarkers.}
\par
The comparative analysis of the heterodimer and FK model in the case of amyloid-$\beta$ shows that using FK provides comparable results in most cases. In particular, the solutions have a slight discrepancy near the average value of the solution, as shown in Figure \ref{fig:sens_beta_FK}. {Indeed, the major importance of the parameter $p_\Delta$ on the final solution, and the sigmoidal shapes of the functions, justifies the FK model in the numerical modeling of the amyloid-$\beta$. This is coherent both with medical literature \cite{van_oostveen_imaging_2021} shapes and with the mathematical works on the topic \cite{weickenmeierPhysicsbasedModelExplains2019,corti_discontinuous_2023}.} These findings suggest a good level of accuracy in describing the phenomenon of the amyloid-$\beta$ spreading independently on the chosen model.
\subsection{Limitations}
\label{sec:limitations}
The main limitation of this work is the two-dimensional domains. The extension of the analysis to three-dimensional complete brain geometries can give new insights into the accuracy of the models, also introducing the heterogeneity in the parameters. Also, the calibration of the model parameters using localized measures of concentrations and not a homogeneous one could be an improvement. However, this would require many \textit{post-mortem} biological analyses of different brain regions or a deeper analysis of correlations between concentration values and PET images. Another limitation is the absence of a sensitivity analysis of $k_{12}$, which we neglected because it is not computable from biological measures. However, we underline that changes in this value are not related to changes in the equilibrium nature but only to the velocity of the pathology development. {Finally, the last limitation is the assumption of the distribution, which would require a detailed inferential estimation of a large number of data. In particular, for what concerns the amyloid-$\beta$ a derivation of the distribution from raw data would be a great improvement for the work.}
\section{Conclusions}
\label{sec:conclusions}
In this work, we performed a sensitivity analysis of heterodimer and Fisher-Kolmogorov models, considering the impact of variations in the protein concentrations at the equilibrium on the problem solution. The parameter values were derived from biological measurements in the brain cortex for the tau protein and the amyloid-$\beta$ in Alzheimer's diseased patients. Finally, we compared the results of the two models, discussing the importance of the model choice, arising, in particular, in the simulation of tau spreading, where some oscillations of the concentration values seem to happen before reaching the equilibrium solution in one model and not in the other. {This work confirms the quality of the approximations of the amyloid-$\beta$ spreading by means of the Fisher-Kolmogorov equation, also due to the major importance of $p_\Delta$ deriving from the sensitivity analysis.} An interesting future development of the work is the study of three-dimensional brain geometries, trying to locally calibrate the parameters using the correlation between positron emission tomography images and protein concentrations.
\section*{Acknowledgments}
The brain MRI images were provided by OASIS-3: Longitudinal Multimodal Neuroimaging: Principal Investigators: T. Benzinger, D. Marcus, J. Morris; NIH P30 AG066444, P50 AG00561, P30 NS09857781, P01 AG026276, P01 AG003991, R01 AG043434, UL1 TR000448, R01 EB009352. AV-45 doses were provided by Avid Radiopharmaceuticals, a wholly-owned subsidiary of Eli Lilly.

\section*{Declaration of competing interests}
The author declares that he has no known competing financial interests or personal relationships that could have appeared to influence the work reported in this article.

\bibliographystyle{hieeetr}
\bibliography{sample.bib}
\end{document}